\documentclass[1pt]{amsart}

\usepackage{amsmath,amsthm, amscd, amssymb, amsfonts}

\usepackage{hhline}

\newcommand{\Ind}{\operatorname{Ind}}
\newcommand{\ord}{\operatorname{ord}}
\newcommand{\oper}{\#}

\DeclareMathOperator*{\prode}{\boxtimes}

\newcommand\toba{{\mathfrak B }}
\newcommand\trasp{\pi}
\newcommand{\gr}{\operatorname{gr}}
\newcommand{\trid}{\triangleright}

\newcommand\compvert[2]{\genfrac{}{}{-1pt}{0}{#1}{#2}}

\newcommand\mvert[2]{\begin{tiny}\begin{matrix}#1\vspace{-4pt}\\#2\end{matrix}\end{tiny}}
\DeclareMathOperator*{\Tim}{\times}
\newcommand{\Times}[2]{\sideset{_#1}{_#2}\Tim}

\newcommand{\cx}{{\daleth}}

\newcommand{\lu}{{\nu}}
\newcommand{\ld}{{\ell}}

\newcommand{\abofe}{{A(\lu, \ld)}}

\newcommand{\bbofe}{{B(\lu, \ld)}}

\newcommand{\proy}{\Pi}

\newcommand{\Lc}{{\mathcal L}}
\newcommand{\Y}{{\mathcal Y}}
\newcommand{\R}{{\mathcal R}}
\newcommand{\W}{{\mathcal W}}
\newcommand{\Zc}{{\mathcal Z}}
\newcommand{\daga}{{\dagger}}

\newcommand{\acts}{{\rightharpoondown}}

\newcommand{\ku}{\mathbb C}
\newcommand{\K}{{\mathcal K}}
\newcommand{\Z}{{\mathbb Z}}
\newcommand{\N}{{\mathbb N}}
\newcommand{\I}{{\mathbb I}}

\newcommand{\G}{{\mathbb G}}
\newcommand{\M}{{\mathcal M}}
\newcommand{\Q}{{\mathcal Q}}
\newcommand{\F}{{\mathcal F}}
\newcommand{\C}{{\mathcal C}}
\newcommand{\D}{{\mathcal D}}
\newcommand{\Ec}{{\mathbf E}}
\newcommand{\Ee}{{\mathcal E}}
\newcommand{\Kc}{{\mathbf K}}

\newcommand{\uno}{{\bf 1}}
\newcommand{\uv}{{\bf 1_v}}
\newcommand{\uh}{{\bf 1_h}}
\newcommand{\m}{\mathcal{M}}
\newcommand{\n}{\mathcal{N}}

\newcommand{\Dc}{{\mathbf D}}
\newcommand{\B}{{\mathcal B}}
\newcommand{\T}{{\mathcal T}}
\newcommand{\Hc}{{\mathcal H}}
\newcommand{\Vc}{{\mathcal V}}
\newcommand{\Pc}{{\mathcal P}}
\newcommand{\Oc}{{\mathcal O}}
\newcommand{\ydh}{{}^H_H\mathcal{YD}}

\newcommand{\Ss}{{\mathcal S}}
\newcommand{\Vect}{\operatorname{Vec}}
\newcommand{\End}{\operatorname{End}}
\newcommand{\Aut}{\operatorname{Aut}}
\newcommand{\Int}{\operatorname{Int}}
\newcommand{\Ext}{\operatorname{Ext}}
\newcommand\tr{\operatorname{tr}}
\newcommand\Rep{\operatorname{Rep}}

\newcommand\card{\operatorname{card}}
\newcommand{\unosigma}{{\bf 1}^{\sigma}}
\newcommand\tot{\operatorname{tot}}
\newcommand\sgn{\operatorname{sgn}}
\newcommand\ad{\operatorname{ad}}
\newcommand\Hom{\operatorname{Hom}}
\newcommand\opext{\operatorname{Opext}}
\newcommand\Tot{\operatorname{Tot}}
\newcommand\Map{\operatorname{Map}}
\newcommand{\fde}{{\triangleright}}
\newcommand{\fiz}{{\triangleleft}}
\newcommand{\wC}{\widehat{C}}
\newcommand{\la}{\langle}
\newcommand{\ra}{\rangle}
\newcommand{\Comod}{\mbox{\rm Comod\,}}
\newcommand{\Mod}{\mbox{\rm Mod\,}}
\newcommand{\Sg}{{\mathfrak S}}
\newcommand{\funcion}{\mathfrak p}
\newcommand{\tauo}{{\widehat\tau}}
\newcommand{\prin}{t}%\partial_0}
\newcommand{\fin}{b}%\partial_1}
\newcommand{\pri}{r}%\epsilon_0}
\newcommand{\fine}{l}%\epsilon_1}
\newcommand\rh{\sim_{H}}
\newcommand\rv{\sim_{V}}
\newcommand\rd{\sim_{D}}

\newcommand\V{\operatorname{Vec}}

\theoremstyle{plain}

\renewcommand{\baselinestretch}{1.2}

\title{The character tables of centralizers in  Weyl Group of $E_8$
V}
\author{ \small Shouchuan Zhang,     Peng Wang,  Jing Cheng, Hui Yang
}
\address{ Mathematics Department of Hunan University,\newline \indent Changsha China,
410082, E-mail:\\ z9491@yahoo.com.cn }
\date{}

\begin{document}
\newtheorem{Proposition}{\quad Proposition}[section]
\newtheorem{Theorem}{\quad Theorem}
\newtheorem{Definition}[Proposition]{\quad Definition}
\newtheorem{Corollary}[Proposition]{\quad Corollary}
\newtheorem{Lemma}[Proposition]{\quad Lemma}
\newtheorem{Example}[Proposition]{\quad Example}

\maketitle \addtocounter{section}{-1}

\numberwithin{equation}{section}

\date{}

%\title{
%On Pointed Hopf Algebras with Weyl Groups}
%\thanks{}
%\author{ Shouchuan Zhang,     Peng Wang, Jing Cheng, Hui Yang  \\
% Mathematics Department of Hunan University,\\Changsha China,
%410082, E-mail:\\ z9491@yahoo.com.cn }
%\date{}
%\maketitle
%\begin{center}
%\par\ \\
%\xiaosanhao{\hei对称群的分歧数据系统和交换群的元素系统分类}\\
%\end{center}

\begin {abstract} To classify the finite dimensional pointed Hopf
algebras with Weyl group $G$ of $E_8$, we obtain  the
representatives of conjugacy classes of $G$ and  all character
tables of centralizers of these representatives by means of software
GAP. In this paper we only list character table 95--112.


keywords: GAP, Hopf algebra, Weyl group, character.
\end {abstract}

\section{Introduction}\label {s0}

This article is to contribute to the classification of
finite-dimensional complex pointed Hopf algebras  with Weyl groups
of $E_6$, $E_7$, $F_4$, $G_2$.
 Many papers are about the classification of finite dimensional
pointed Hopf algebras, for example,  \cite{AS98b, AS02, AS00, AS05,
He06, AHS08, AG03, AFZ08, AZ07,Gr00,  Fa07, AF06, AF07, ZZC04, ZCZ08}.

In these research  ones need  the centralizers and character tables
of groups. In this paper we obtain   the representatives of
conjugacy classes of Weyl groups of $E_8$  and all character tables
of centralizers of these representatives by means of software GAP.
In this paper we only list character table  95--112.

By the Cartan-Killing classification of simple Lie algebras over complex field the Weyl groups
to
be considered are $W(A_l), (l\geq 1); $
 $W(B_l), (l \geq 2);$
 $W(C_l), (l \geq 2); $
 $W(D_l), (l\geq 4); $
 $W(E_l), (8 \geq l \geq 6); $
 $W(F_l), ( l=4 );$
$W(G_l), (l=2)$.

It is otherwise desirable to do this in view of the importance of Weyl groups in the theories of
Lie groups, Lie algebras and algebraic groups. For example, the irreducible representations of
Weyl groups were obtained by Frobenius, Schur and Young. The conjugace classes of  $W(F_4)$
were obtained by Wall
 \cite {Wa63} and its character tables were obtained by Kondo \cite {Ko65}.
The conjugace classes  and  character tables of $W(E_6),$ $W(E_7)$ and $W(E_8)$ were obtained
by  Frame \cite {Fr51}.
Carter gave a unified description of the conjugace classes of
 Weyl groups of simple Lie algebras  \cite  {Ca72}.

\section {Program}

 By using the following
program in GAP,  we  obtain  the representatives of conjugacy
classes of Weyl groups of $E_6$  and all character tables of
centralizers of these representatives.

gap$>$ L:=SimpleLieAlgebra("E",6,Rationals);;

gap$>$ R:=RootSystem(L);;

 gap$>$ W:=WeylGroup(R);Display(Order(W));

gap $>$ ccl:=ConjugacyClasses(W);;

gap$>$ q:=NrConjugacyClasses(W);; Display (q);

gap$>$ for i in [1..q] do

$>$ r:=Order(Representative(ccl[i]));Display(r);;

$>$ od; gap

$>$ s1:=Representative(ccl[1]);;cen1:=Centralizer(W,s1);;

gap$>$ cl1:=ConjugacyClasses(cen1);

gap$>$ s1:=Representative(ccl[2]);;cen1:=Centralizer(W,s1);;

 gap$>$
cl2:=ConjugacyClasses(cen1);

gap$>$ s1:=Representative(ccl[3]);;cen1:=Centralizer(W,s1);;

gap$>$ cl3:=ConjugacyClasses(cen1);

 gap$>$
s1:=Representative(ccl[4]);;cen1:=Centralizer(W,s1);;

 gap$>$
cl4:=ConjugacyClasses(cen1);

 gap$>$
s1:=Representative(ccl[5]);;cen1:=Centralizer(W,s1);;

gap$>$ cl5:=ConjugacyClasses(cen1);

 gap$>$
s1:=Representative(ccl[6]);;cen1:=Centralizer(W,s1);;

gap$>$ cl6:=ConjugacyClasses(cen1);

gap$>$ s1:=Representative(ccl[7]);;cen1:=Centralizer(W,s1);;

gap$>$ cl7:=ConjugacyClasses(cen1);

gap$>$ s1:=Representative(ccl[8]);;cen1:=Centralizer(W,s1);;

gap$>$ cl8:=ConjugacyClasses(cen1);

gap$>$ s1:=Representative(ccl[9]);;cen1:=Centralizer(W,s1);;

gap$>$ cl9:=ConjugacyClasses(cen1);

gap$>$ s1:=Representative(ccl[10]);;cen1:=Centralizer(W,s1);;

 gap$>$
cl10:=ConjugacyClasses(cen1);

gap$>$ s1:=Representative(ccl[11]);;cen1:=Centralizer(W,s1);;

gap$>$ cl11:=ConjugacyClasses(cen1);

gap$>$ s1:=Representative(ccl[12]);;cen1:=Centralizer(W,s1);;

gap$>$ cl2:=ConjugacyClasses(cen1);

gap$>$ s1:=Representative(ccl[13]);;cen1:=Centralizer(W,s1);;

 gap$>$
cl13:=ConjugacyClasses(cen1);

 gap$>$
s1:=Representative(ccl[14]);;cen1:=Centralizer(W,s1);;

gap$>$ cl14:=ConjugacyClasses(cen1);

gap$>$ s1:=Representative(ccl[15]);;cen1:=Centralizer(W,s1);;

gap$>$ cl15:=ConjugacyClasses(cen1);

gap$>$ s1:=Representative(ccl[16]);;cen1:=Centralizer(W,s1);;

gap$>$ cl16:=ConjugacyClasses(cen1);

gap$>$ s1:=Representative(ccl[17]);;cen1:=Centralizer(W,s1);;

gap$>$ cl17:=ConjugacyClasses(cen1);

gap$>$ s1:=Representative(ccl[18]);;cen1:=Centralizer(W,s1);;

gap$>$ cl18:=ConjugacyClasses(cen1);

gap$>$ s1:=Representative(ccl[19]);;cen1:=Centralizer(W,s1);;

gap$>$ cl19:=ConjugacyClasses(cen1);

gap$>$ s1:=Representative(ccl[20]);;cen1:=Centralizer(W,s1);;

gap$>$ cl20:=ConjugacyClasses(cen1);

gap$>$ s1:=Representative(ccl[21]);;cen1:=Centralizer(W,s1);;

 gap$>$
cl21:=ConjugacyClasses(cen1);

gap$>$ s1:=Representative(ccl[22]);;cen1:=Centralizer(W,s1);;

 gap$>$
cl22:=ConjugacyClasses(cen1);

gap$>$ s1:=Representative(ccl[23]);;cen1:=Centralizer(W,s1);;

gap$>$ cl23:=ConjugacyClasses(cen1);

gap$>$ s1:=Representative(ccl[24]);;cen1:=Centralizer(W,s1);;

$>$ cl24:=ConjugacyClasses(cen1);

gap$>$ s1:=Representative(ccl[25]);;cen1:=Centralizer(W,s1);

gap$>$ cl25:=ConjugacyClasses(cen1);

gap$>$ for i in [1..q] do

$>$ s:=Representative(ccl[i]);;cen:=Centralizer(W,s);;

$>$ char:=CharacterTable(cen);;Display (cen);Display(char);

 $>$ od;

The programs for Weyl groups of $E_7$, $E_8$, $F_4$ and $ G_2$ are
similar.

{\small

\section {$E_8$}
%%\end {document}
The generators of $G^{s_{95}}$  are:\\
$\left(
\begin{array}{cccccccc} -1&1&1&-2& 1&0&1&0 \\
0&2&0&-3&2&0&1&0 \\ -1&2&1&-4&3&0&1&
0 \\
        -1&4&1&-6&4&0&1&1 \\ -1&3&1&-5&4&-1&1&1 \\ -1&2&1&-4&3&0&0&1 \\
        0&1&1&-3&2&0&0&1 \\ 0&1&1&-2&1&0&0&0 \end{array}\right)$,
  $\left(\begin{array}{cccccccc} 1&1&-1&-1&0&2&-1&0 \\ 0&2&-1&-2&1&2&-1&0 \\ 0&2&-1&-2&0&4&-2&0 \\
        0&3&-2&-3&1&5&-3&1 \\ 0&2&-2&-2&1&4&-3&1 \\ 0&1&-1&-2&1&3&-2&1 \\
        0&1&-1&-1&0&2&-1&1 \\ 0&1&0&-1&0&1&-1&1 \end{array}\right)$.

The representatives of conjugacy classes of   $G^{s_{95}}$ are:\\
$\left(\begin{array}{cccccccc} 1&0&0&0&0&0&0&0 \\ 0&1&0&0&
0&0&0&0 \\ 0&0&1&0&0&0&0&0 \\
        0&0&0&1&0&0&0&0 \\ 0&0&0&0&1&0&0&0 \\ 0&0&0&0&0&1&0&0 \\
        0&0&0&0&0&0&1&0 \\ 0&0&0&0&0&0&0&1 \end{array}\right)$,
  $\left(\begin{array}{cccccccc} -1&-1&1&1&0&-2&1&0 \\ -1&-2&1&2&0&-2&0&0 \\ -1&-3&2&2&0&-3&1&0 \\
        -2&-4&3&3&0&-4&0&1 \\ -2&-3&3&2&0&-3&0&1 \\
        -1&-2&2&2&-1&-2&0&1 \\ 0&-1&1&1&0&-2&0&1 \\ 0&0&1&0&0&-1&0&0 \end{array}\right)$,
 $\left(\begin{array}{cccccccc} -1&-1&1&1&0&-2&1&0 \\ -1&-2&1&2&0&-2&0&0 \\
        -1&-3&2&2&0&-3&1&0 \\ -2&-4&3&3&0&-4&1&-1 \\
        -2&-3&3&2&0&-3&1&-1 \\ -1&-2&2&2&-1&-2&1&-1 \\
        0&-1&1&1&0&-2&1&-1 \\ 0&-1&0&1&0&-1&1&-1
        \end{array}\right)$,  \\
  $\left(\begin{array}{cccccccc} -1&-1&1&1&0&-2&1&0 \\ 0&-2&1&2&-1&-2&1&0 \\ 0&-2&1&2&0&-4&2&0 \\
        0&-4&1&4&-1&-5&3&0 \\ 0&-3&1&3&-1&-4&3&0 \\ 0&-2&0&3&-1&-3&2&0 \\
        0&-2&0&2&0&-2&1&0 \\ 0&-1&0&1&0&-1&0&1 \end{array}\right)$,
  $\left(\begin{array}{cccccccc} -1&-1&1&1&0&-2&1&0 \\ 0&-2&1&2&-1&-2&1&0 \\ 0&-2&1&2&0&-4&2&0 \\
        0&-3&2&3&-1&-5&2&1 \\ 0&-2&2&2&-1&-4&2&1 \\ 0&-1&1&2&-1&-3&1&1 \\
        0&-1&1&1&0&-2&0&1 \\ 0&0&1&0&0&-1&0&0 \end{array}\right)$,
  $\left(\begin{array}{cccccccc} -1&-1&1&1&0&-2&1&0 \\ 0&-2&1&2&-1&-2&1&0 \\ 0&-2&1&2&0&-4&2&0 \\
        0&-3&2&3&-1&-5&3&-1 \\ 0&-2&2&2&-1&-4&3&-1 \\
        0&-1&1&2&-1&-3&2&-1 \\ 0&-1&1&1&0&-2&1&-1 \\
        0&-1&0&1&0&-1&1&-1 \end{array}\right)$,  \\
  $\left(\begin{array}{cccccccc} -1&-1&1&1&0&-2&1&0 \\ 0&-1&0&2&0&-3&1&0 \\ -1&-2&1&2&1&-4&1&0 \\
        -1&-3&0&4&1&-6&2&0 \\ -1&-2&0&3&1&-5&2&0 \\ -1&-2&0&3&0&-3&1&0 \\
        0&-2&0&2&0&-2&1&0 \\ 0&-1&0&1&0&-1&0&1 \end{array}\right)$,
  $\left(\begin{array}{cccccccc} -1&-1&1&1&0&-2&1&0 \\ 0&-1&0&2&0&-3&1&0 \\ -1&-2&1&2&1&-4&1&0 \\
        -1&-2&1&3&1&-6&1&1 \\ -1&-1&1&2&1&-5&1&1 \\ -1&-1&1&2&0&-3&0&1 \\
        0&-1&1&1&0&-2&0&1 \\ 0&0&1&0&0&-1&0&0 \end{array}\right)$,
  $\left(\begin{array}{cccccccc} -1&-1&1&1&0&-2&1&0 \\ 0&-1&0&2&0&-3&1&0 \\ -1&-2&1&2&1&-4&1&0 \\
        -1&-2&1&3&1&-6&2&-1 \\ -1&-1&1&2&1&-5&2&-1 \\
        -1&-1&1&2&0&-3&1&-1 \\ 0&-1&1&1&0&-2&1&-1 \\
        0&-1&0&1&0&-1&1&-1 \end{array}\right)$,  \\
  $\left(\begin{array}{cccccccc} -1&0&0&0&0&0&0&0 \\ -1&-1&0&0&1&0&-1&0 \\ -1&-1&0&0&0&1&-1&0 \\
        -2&-2&0&0&1&1&-2&1 \\ -2&-2&0&1&0&1&-2&1 \\ -1&-2&0&1&0&0&-1&1 \\
        0&-1&-1&1&0&0&-1&1 \\ 0&0&0&0&0&0&-1&1 \end{array}\right)$,
  $\left(\begin{array}{cccccccc} -1&0&0&0&0&0&0&0 \\ -1&-1&0&0&1&0&-1&0 \\ -1&-1&0&0&0&1&-1&0 \\
        -2&-2&0&0&1&1&-1&-1 \\ -2&-2&0&1&0&1&-1&-1 \\
        -1&-2&0&1&0&0&0&-1 \\ 0&-1&-1&1&0&0&0&-1 \\ 0&-1&-1&1&0&0&0&0 \end{array}\right)$,
  $\left(\begin{array}{cccccccc} -1&0&0&0&0&0&0&0 \\ -1&-1&0&0&1&0&-1&0 \\ -1&-1&0&0&0&1&-1&0 \\
        -2&-1&1&-1&1&1&-2&0 \\ -2&-1&1&0&0&1&-2&0 \\
        -1&-1&1&0&0&0&-1&0 \\ 0&0&0&0&0&0&-1&0 \\ 0&0&0&0&0&0&0&-1
        \end{array}\right)$,  \\
  $\left(\begin{array}{cccccccc} -1&0&0&0&0&0&0&0 \\ 0&-1&0&0&0&0&0&0 \\ 0&0&-1&0&0&0&0&0 \\
        0&-1&-1&0&0&0&0&1 \\ 0&-1&-1&1&-1&0&0&1 \\ 0&-1&-1&1&0&-1&0&1 \\
        0&-1&-1&1&0&0&-1&1 \\ 0&0&0&0&0&0&-1&1 \end{array}\right)$,
  $\left(\begin{array}{cccccccc} -1&0&0&0&0&0&0&0 \\ 0&-1&0&0&0&0&0&0 \\ 0&0&-1&0&0&0&0&0 \\
        0&-1&-1&0&0&0&1&-1 \\ 0&-1&-1&1&-1&0&1&-1 \\
        0&-1&-1&1&0&-1&1&-1 \\ 0&-1&-1&1&0&0&0&-1 \\ 0&-1&-1&1&0&0&0&0 \end{array}\right)$,
  $\left(\begin{array}{cccccccc} -1&0&0&0&0&0&0&0 \\ 0&-1&0&0&0&0&0&0 \\ 0&0&-1&0&0&0&0&0 \\
        0&0&0&-1&0&0&0&0 \\ 0&0&0&0&-1&0&0&0 \\ 0&0&0&0&0&-1&0&0 \\
        0&0&0&0&0&0&-1&0 \\ 0&0&0&0&0&0&0&-1 \end{array}\right)$,
        \\
  $\left(\begin{array}{cccccccc} -1&0&0&0&0&0&0&0 \\ 0&0&-1&0&1&-1&0&0 \\ -1&0&-1&0&1&0&-1&0 \\
        -1&0&-2&0&2&-1&-1&1 \\ -1&0&-2&1&1&-1&-1&1 \\
        -1&-1&-1&1&1&-1&-1&1 \\ 0&-1&-1&1&0&0&-1&1 \\ 0&0&0&0&0&0&-1&1 \end{array}\right)$,
  $\left(\begin{array}{cccccccc} -1&0&0&0&0&0&0&0 \\ 0&0&-1&0&1&-1&0&0 \\ -1&0&-1&0&1&0&-1&0 \\
        -1&0&-2&0&2&-1&0&-1 \\ -1&0&-2&1&1&-1&0&-1 \\
        -1&-1&-1&1&1&-1&0&-1 \\ 0&-1&-1&1&0&0&0&-1 \\
        0&-1&-1&1&0&0&0&0 \end{array}\right)$,
  $\left(\begin{array}{cccccccc} -1&0&0&0&0&0&0&0 \\ 0&0&-1&0&1&-1&0&0 \\ -1&0&-1&0&1&0&-1&0 \\
        -1&1&-1&-1&2&-1&-1&0 \\ -1&1&-1&0&1&-1&-1&0 \\
        -1&0&0&0&1&-1&-1&0 \\ 0&0&0&0&0&0&-1&0 \\ 0&0&0&0&0&0&0&-1
        \end{array}\right)$,  \\
  $\left(\begin{array}{cccccccc} -1&1&0&1&-2&1&0&0 \\ -1&1&0&1&-2&2&-1&0 \\ -1&1&0&2&-4&3&-1&0 \\
        -2&1&0&3&-5&4&-2&1 \\ -2&1&0&2&-4&4&-2&1 \\ -1&1&0&1&-3&3&-1&1 \\
        0&1&-1&1&-2&2&-1&1 \\ 0&1&0&0&-1&1&-1&1 \end{array}\right)$,
  $\left(\begin{array}{cccccccc} -1&1&0&1&-2&1&0&0 \\ -1&1&0&1&-2&2&-1&0 \\ -1&1&0&2&-4&3&-1&0 \\
        -2&1&0&3&-5&4&-1&-1 \\ -2&1&0&2&-4&4&-1&-1 \\
        -1&1&0&1&-3&3&0&-1 \\ 0&1&-1&1&-2&2&0&-1 \\ 0&0&-1&1&-1&1&0&0 \end{array}\right)$,
  $\left(\begin{array}{cccccccc} -1&1&0&1&-2&1&0&0 \\ -1&1&0&1&-2&2&-1&0 \\ -1&1&0&2&-4&3&-1&0 \\
        -2&2&1&2&-5&4&-2&0 \\ -2&2&1&1&-4&4&-2&0 \\ -1&2&1&0&-3&3&-1&0 \\
        0&2&0&0&-2&2&-1&0 \\ 0&1&0&0&-1&1&0&-1 \end{array}\right)$,
        \\
  $\left(\begin{array}{cccccccc} -1&1&0&1&-2&1&0&0 \\ 0&1&0&1&-3&2&0&0 \\ 0&2&-1&2&-4&2&0&0 \\
        0&2&-1&3&-6&3&0&1 \\ 0&2&-1&2&-5&3&0&1 \\ 0&2&-1&1&-3&2&0&1 \\
        0&1&-1&1&-2&2&-1&1 \\ 0&1&0&0&-1&1&-1&1 \end{array}\right)$,
  $\left(\begin{array}{cccccccc} -1&1&0&1&-2&1&0&0 \\ 0&1&0&1&-3&2&0&0 \\ 0&2&-1&2&-4&2&0&0 \\
        0&2&-1&3&-6&3&1&-1 \\ 0&2&-1&2&-5&3&1&-1 \\ 0&2&-1&1&-3&2&1&-1 \\
        0&1&-1&1&-2&2&0&-1 \\ 0&0&-1&1&-1&1&0&0 \end{array}\right)$,
  $\left(\begin{array}{cccccccc} -1&1&0&1&-2&1&0&0 \\ 0&1&0&1&-3&2&0&0 \\ 0&2&-1&2&-4&2&0&0 \\
        0&3&0&2&-6&3&0&0 \\ 0&3&0&1&-5&3&0&0 \\ 0&3&0&0&-3&2&0&0 \\
        0&2&0&0&-2&2&-1&0 \\ 0&1&0&0&-1&1&0&-1 \end{array}\right)$,
        \\
  $\left(\begin{array}{cccccccc} -1&1&0&1&-2&1&0&0 \\ 0&2&-1&1&-2&1&0&0 \\ -1&2&-1&2&-3&2&-1&0 \\
        -1&3&-2&3&-4&2&-1&1 \\ -1&3&-2&2&-3&2&-1&1 \\
        -1&2&-1&1&-2&2&-1&1 \\ 0&1&-1&1&-2&2&-1&1 \\ 0&1&0&0&-1&1&-1&1 \end{array}\right)$,
  $\left(\begin{array}{cccccccc} -1&1&0&1&-2&1&0&0 \\ 0&2&-1&1&-2&1&0&0 \\ -1&2&-1&2&-3&2&-1&0 \\
        -1&3&-2&3&-4&2&0&-1 \\ -1&3&-2&2&-3&2&0&-1 \\
        -1&2&-1&1&-2&2&0&-1 \\ 0&1&-1&1&-2&2&0&-1 \\ 0&0&-1&1&-1&1&0&0 \end{array}\right)$,
  $\left(\begin{array}{cccccccc} -1&1&0&1&-2&1&0&0 \\ 0&2&-1&1&-2&1&0&0 \\ -1&2&-1&2&-3&2&-1&0 \\
        -1&4&-1&2&-4&2&-1&0 \\ -1&4&-1&1&-3&2&-1&0 \\
        -1&3&0&0&-2&2&-1&0 \\ 0&2&0&0&-2&2&-1&0 \\ 0&1&0&0&-1&1&0&-1
        \end{array}\right)$,  \\
  $\left(\begin{array}{cccccccc} -1&1&1&-2&1&0&1&0 \\ -1&1&1&-3&2&1&0&0 \\ -1&1&2&-4&2&1&1&0 \\
        -2&1&2&-5&3&2&1&0 \\ -2&0&2&-4&3&1&1&0 \\ -1&0&1&-3&2&1&1&0 \\
        0&0&0&-2&2&0&1&0 \\ 0&0&0&-1&1&0&0&1 \end{array}\right)$,
  $\left(\begin{array}{cccccccc} -1&1&1&-2&1&0&1&0 \\ -1&1&1&-3&2&1&0&0 \\ -1&1&2&-4&2&1&1&0 \\
        -2&2&3&-6&3&2&0&1 \\ -2&1&3&-5&3&1&0&1 \\ -1&1&2&-4&2&1&0&1 \\
        0&1&1&-3&2&0&0&1 \\ 0&1&1&-2&1&0&0&0 \end{array}\right)$,
  $\left(\begin{array}{cccccccc} -1&1&1&-2&1&0&1&0 \\ -1&1&1&-3&2&1&0&0 \\ -1&1&2&-4&2&1&1&0 \\
        -2&2&3&-6&3&2&1&-1 \\ -2&1&3&-5&3&1&1&-1 \\ -1&1&2&-4&2&1&1&-1 \\
        0&1&1&-3&2&0&1&-1 \\ 0&0&0&-1&1&0&1&-1 \end{array}\right)$,
        \\
  $\left(\begin{array}{cccccccc} -1&1&1&-2&1&0&1&0 \\ 0&1&1&-3&1&1&1&0 \\ 0&2&1&-4&2&0&2&0 \\
        0&2&1&-5&2&1&3&0 \\ 0&1&1&-4&2&0&3&0 \\ 0&1&0&-3&2&0&2&0 \\
        0&0&0&-2&2&0&1&0 \\ 0&0&0&-1&1&0&0&1 \end{array}\right)$,
  $\left(\begin{array}{cccccccc} -1&1&1&-2&1&0&1&0 \\ 0&1&1&-3&1&1&1&0 \\ 0&2&1&-4&2&0&2&0 \\
        0&3&2&-6&2&1&2&1 \\ 0&2&2&-5&2&0&2&1 \\ 0&2&1&-4&2&0&1&1 \\
        0&1&1&-3&2&0&0&1 \\ 0&1&1&-2&1&0&0&0 \end{array}\right)$,
  $\left(\begin{array}{cccccccc} -1&1&1&-2&1&0&1&0 \\ 0&1&1&-3&1&1&1&0 \\ 0&2&1&-4&2&0&2&0 \\
        0&3&2&-6&2&1&3&-1 \\ 0&2&2&-5&2&0&3&-1 \\ 0&2&1&-4&2&0&2&-1 \\
        0&1&1&-3&2&0&1&-1 \\ 0&0&0&-1&1&0&1&-1 \end{array}\right)$,
        \\
  $\left(\begin{array}{cccccccc} -1&1&1&-2&1&0&1&0 \\ 0&2&0&-3&2&0&1&0 \\ -1&2&1&-4&3&0&1&0 \\
        -1&3&0&-5&4&0&2&0 \\ -1&2&0&-4&4&-1&2&0 \\ -1&1&0&-3&3&0&1&0 \\
        0&0&0&-2&2&0&1&0 \\ 0&0&0&-1&1&0&0&1 \end{array}\right)$,
  $\left(\begin{array}{cccccccc} -1&1&1&-2&1&0&1&0 \\ 0&2&0&-3&2&0&1&0 \\ -1&2&1&-4&3&0&1&0 \\
        -1&4&1&-6&4&0&1&1 \\ -1&3&1&-5&4&-1&1&1 \\ -1&2&1&-4&3&0&0&1 \\
        0&1&1&-3&2&0&0&1 \\ 0&1&1&-2&1&0&0&0 \end{array}\right)$,
  $\left(\begin{array}{cccccccc} -1&1&1&-2&1&0&1&0 \\ 0&2&0&-3&2&0&1&0 \\ -1&2&1&-4&3&0&1&0 \\
        -1&4&1&-6&4&0&2&-1 \\ -1&3&1&-5&4&-1&2&-1 \\
        -1&2&1&-4&3&0&1&-1 \\ 0&1&1&-3&2&0&1&-1 \\ 0&0&0&-1&1&0&1&-1
        \end{array}\right)$,  \\
   $\left(\begin{array}{cccccccc} 1&-1&-1&2&-1&0&-1&0 \\ 0&-2&0&3&-2&0&-1&0 \\
        1&-2&-1&4&-3&0&-1&0 \\ 1&-4&-1&6&-4&0&-2&1 \\
        1&-3&-1&5&-4&1&-2&1 \\ 1&-2&-1&4&-3&0&-1&1 \\
        0&-1&-1&3&-2&0&-1&1 \\ 0&0&0&1&-1&0&-1&1 \end{array}\right)$,
  $\left(\begin{array}{cccccccc} 1&-1&-1&2&-1&0&-1&0 \\ 0&-2&0&3&-2&0&-1&0 \\ 1&-2&-1&4&-3&0&-1&0 \\
        1&-4&-1&6&-4&0&-1&-1 \\ 1&-3&-1&5&-4&1&-1&-1 \\
        1&-2&-1&4&-3&0&0&-1 \\ 0&-1&-1&3&-2&0&0&-1 \\
        0&-1&-1&2&-1&0&0&0 \end{array}\right)$,
  $\left(\begin{array}{cccccccc} 1&-1&-1&2&-1&0&-1&0 \\ 0&-2&0&3&-2&0&-1&0 \\ 1&-2&-1&4&-3&0&-1&0 \\
        1&-3&0&5&-4&0&-2&0 \\ 1&-2&0&4&-4&1&-2&0 \\ 1&-1&0&3&-3&0&-1&0 \\
        0&0&0&2&-2&0&-1&0 \\ 0&0&0&1&-1&0&0&-1 \end{array}\right)$,
        \\
  $\left(\begin{array}{cccccccc} 1&-1&-1&2&-1&0&-1&0 \\ 0&-1&-1&3&-1&-1&-1&0 \\
        0&-2&-1&4&-2&0&-2&0 \\ 0&-3&-2&6&-2&-1&-3&1 \\
        0&-2&-2&5&-2&0&-3&1 \\ 0&-2&-1&4&-2&0&-2&1 \\
        0&-1&-1&3&-2&0&-1&1 \\ 0&0&0&1&-1&0&-1&1 \end{array}\right)$,
  $\left(\begin{array}{cccccccc} 1&-1&-1&2&-1&0&-1&0 \\ 0&-1&-1&3&-1&-1&-1&0 \\
        0&-2&-1&4&-2&0&-2&0 \\ 0&-3&-2&6&-2&-1&-2&-1 \\
        0&-2&-2&5&-2&0&-2&-1 \\ 0&-2&-1&4&-2&0&-1&-1 \\
        0&-1&-1&3&-2&0&0&-1 \\ 0&-1&-1&2&-1&0&0&0 \end{array}\right)$,
  $\left(\begin{array}{cccccccc} 1&-1&-1&2&-1&0&-1&0 \\ 0&-1&-1&3&-1&-1&-1&0 \\
        0&-2&-1&4&-2&0&-2&0 \\ 0&-2&-1&5&-2&-1&-3&0 \\
        0&-1&-1&4&-2&0&-3&0 \\ 0&-1&0&3&-2&0&-2&0 \\ 0&0&0&2&-2&0&-1&0 \\
        0&0&0&1&-1&0&0&-1 \end{array}\right)$,  \\
  $\left(\begin{array}{cccccccc} 1&-1&-1&2&-1&0&-1&0 \\ 1&-1&-1&3&-2&-1&0&0 \\
        1&-1&-2&4&-2&-1&-1&0 \\ 2&-2&-3&6&-3&-2&-1&1 \\
        2&-1&-3&5&-3&-1&-1&1 \\ 1&-1&-2&4&-2&-1&-1&1 \\
        0&-1&-1&3&-2&0&-1&1 \\ 0&0&0&1&-1&0&-1&1 \end{array}\right)$,
  $\left(\begin{array}{cccccccc} 1&-1&-1&2&-1&0&-1&0 \\ 1&-1&-1&3&-2&-1&0&0 \\
        1&-1&-2&4&-2&-1&-1&0 \\ 2&-2&-3&6&-3&-2&0&-1 \\
        2&-1&-3&5&-3&-1&0&-1 \\ 1&-1&-2&4&-2&-1&0&-1 \\
        0&-1&-1&3&-2&0&0&-1 \\ 0&-1&-1&2&-1&0&0&0 \end{array}\right)$,
  $\left(\begin{array}{cccccccc} 1&-1&-1&2&-1&0&-1&0 \\ 1&-1&-1&3&-2&-1&0&0 \\
        1&-1&-2&4&-2&-1&-1&0 \\ 2&-1&-2&5&-3&-2&-1&0 \\
        2&0&-2&4&-3&-1&-1&0 \\ 1&0&-1&3&-2&-1&-1&0 \\
        0&0&0&2&-2&0&-1&0 \\ 0&0&0&1&-1&0&0&-1 \end{array}\right)$,
        \\
  $\left(\begin{array}{cccccccc} 1&-1&0&-1&2&-1&0&0 \\ 0&-2&1&-1&2&-1&0&0 \\ 1&-2&1&-2&3&-2&1&0 \\
        1&-4&1&-2&4&-2&1&0 \\ 1&-4&1&-1&3&-2&1&0 \\ 1&-3&0&0&2&-2&1&0 \\
        0&-2&0&0&2&-2&1&0 \\ 0&-1&0&0&1&-1&0&1 \end{array}\right)$,
  $\left(\begin{array}{cccccccc} 1&-1&0&-1&2&-1&0&0 \\ 0&-2&1&-1&2&-1&0&0 \\ 1&-2&1&-2&3&-2&1&0 \\
        1&-3&2&-3&4&-2&0&1 \\ 1&-3&2&-2&3&-2&0&1 \\ 1&-2&1&-1&2&-2&0&1 \\
        0&-1&1&-1&2&-2&0&1 \\ 0&0&1&-1&1&-1&0&0 \end{array}\right)$,
  $\left(\begin{array}{cccccccc} 1&-1&0&-1&2&-1&0&0 \\ 0&-2&1&-1&2&-1&0&0 \\ 1&-2&1&-2&3&-2&1&0 \\
        1&-3&2&-3&4&-2&1&-1 \\ 1&-3&2&-2&3&-2&1&-1 \\
        1&-2&1&-1&2&-2&1&-1 \\ 0&-1&1&-1&2&-2&1&-1 \\
        0&-1&0&0&1&-1&1&-1 \end{array}\right)$,  \\
  $\left(\begin{array}{cccccccc} 1&-1&0&-1&2&-1&0&0 \\ 0&-1&0&-1&3&-2&0&0 \\ 0&-2&1&-2&4&-2&0&0 \\
        0&-3&0&-2&6&-3&0&0 \\ 0&-3&0&-1&5&-3&0&0 \\ 0&-3&0&0&3&-2&0&0 \\
        0&-2&0&0&2&-2&1&0 \\ 0&-1&0&0&1&-1&0&1 \end{array}\right)$,
  $\left(\begin{array}{cccccccc} 1&-1&0&-1&2&-1&0&0 \\ 0&-1&0&-1&3&-2&0&0 \\ 0&-2&1&-2&4&-2&0&0 \\
        0&-2&1&-3&6&-3&-1&1 \\ 0&-2&1&-2&5&-3&-1&1 \\
        0&-2&1&-1&3&-2&-1&1 \\ 0&-1&1&-1&2&-2&0&1 \\ 0&0&1&-1&1&-1&0&0 \end{array}\right)$,
  $\left(\begin{array}{cccccccc} 1&-1&0&-1&2&-1&0&0 \\ 0&-1&0&-1&3&-2&0&0 \\
        0&-2&1&-2&4&-2&0&0 \\ 0&-2&1&-3&6&-3&0&-1 \\
        0&-2&1&-2&5&-3&0&-1 \\ 0&-2&1&-1&3&-2&0&-1 \\
        0&-1&1&-1&2&-2&1&-1 \\ 0&-1&0&0&1&-1&1&-1
        \end{array}\right)$,  \\
  $\left(\begin{array}{cccccccc} 1&-1&0&-1&2&-1&0&0 \\ 1&-1&0&-1&2&-2&1&0 \\ 1&-1&0&-2&4&-3&1&0 \\
        2&-2&-1&-2&5&-4&2&0 \\ 2&-2&-1&-1&4&-4&2&0 \\
        1&-2&-1&0&3&-3&1&0 \\ 0&-2&0&0&2&-2&1&0 \\ 0&-1&0&0&1&-1&0&1 \end{array}\right)$,
  $\left(\begin{array}{cccccccc} 1&-1&0&-1&2&-1&0&0 \\ 1&-1&0&-1&2&-2&1&0 \\
        1&-1&0&-2&4&-3&1&0 \\ 2&-1&0&-3&5&-4&1&1 \\ 2&-1&0&-2&4&-4&1&1 \\
        1&-1&0&-1&3&-3&0&1 \\ 0&-1&1&-1&2&-2&0&1 \\ 0&0&1&-1&1&-1&0&0 \end{array}\right)$,
  $\left(\begin{array}{cccccccc} 1&-1&0&-1&2&-1&0&0 \\ 1&-1&0&-1&2&-2&1&0 \\
        1&-1&0&-2&4&-3&1&0 \\ 2&-1&0&-3&5&-4&2&-1 \\
        2&-1&0&-2&4&-4&2&-1 \\ 1&-1&0&-1&3&-3&1&-1 \\
        0&-1&1&-1&2&-2&1&-1 \\ 0&-1&0&0&1&-1&1&-1
        \end{array}\right)$,  \\
  $\left(\begin{array}{cccccccc} 1&0&0&0&0&0&0&0 \\ 0&0&1&0&-1&1&0&0 \\ 1&0&1&0&-1&0&1&0 \\
        1&-1&1&1&-2&1&1&0 \\ 1&-1&1&0&-1&1&1&0 \\ 1&0&0&0&-1&1&1&0 \\
        0&0&0&0&0&0&1&0 \\ 0&0&0&0&0&0&0&1 \end{array}\right)$,
  $\left(\begin{array}{cccccccc} 1&0&0&0&0&0&0&0 \\ 0&0&1&0&-1&1&0&0 \\ 1&0&1&0&-1&0&1&0 \\
        1&0&2&0&-2&1&0&1 \\ 1&0&2&-1&-1&1&0&1 \\ 1&1&1&-1&-1&1&0&1 \\
        0&1&1&-1&0&0&0&1 \\ 0&1&1&-1&0&0&0&0 \end{array}\right)$,
  $\left(\begin{array}{cccccccc} 1&0&0&0&0&0&0&0 \\ 0&0&1&0&-1&1&0&0 \\ 1&0&1&0&-1&0&1&0 \\
        1&0&2&0&-2&1&1&-1 \\ 1&0&2&-1&-1&1&1&-1 \\ 1&1&1&-1&-1&1&1&-1 \\
        0&1&1&-1&0&0&1&-1 \\ 0&0&0&0&0&0&1&-1 \end{array}\right)$,
        \\
  $\left(\begin{array}{cccccccc} -1&-1&1&1&0&-2&1&0 \\ -1&-2&1&2&0&-2&0&0 \\ -1&-3&2&2&0&-3&1&0 \\
        -2&-5&2&4&0&-4&1&0 \\ -2&-4&2&3&0&-3&1&0 \\
        -1&-3&1&3&-1&-2&1&0 \\ 0&-2&0&2&0&-2&1&0 \\ 0&-1&0&1&0&-1&0&1 \end{array}\right)$,
  $\left(\begin{array}{cccccccc} 1&0&0&0&0&0&0&0 \\ 0&1&0&0&0&0&0&0 \\ 0&0&1&0&0&0&0&0 \\
        0&1&1&0&0&0&-1&1 \\ 0&1&1&-1&1&0&-1&1 \\ 0&1&1&-1&0&1&-1&1 \\
        0&1&1&-1&0&0&0&1 \\ 0&1&1&-1&0&0&0&0 \end{array}\right)$,
  $\left(\begin{array}{cccccccc} 1&0&0&0&0&0&0&0 \\ 0&1&0&0&0&0&0&0 \\ 0&0&1&0&0&0&0&0 \\
        0&1&1&0&0&0&0&-1 \\ 0&1&1&-1&1&0&0&-1 \\ 0&1&1&-1&0&1&0&-1 \\
        0&1&1&-1&0&0&1&-1 \\ 0&0&0&0&0&0&1&-1 \end{array}\right)$,
        \\
  $\left(\begin{array}{cccccccc} 1&0&0&0&0&0&0&0 \\ 1&1&0&0&-1&0&1&0 \\ 1&1&0&0&0&-1&1&0 \\
        2&1&-1&1&-1&-1&2&0 \\ 2&1&-1&0&0&-1&2&0 \\ 1&1&-1&0&0&0&1&0 \\
        0&0&0&0&0&0&1&0 \\ 0&0&0&0&0&0&0&1 \end{array}\right)$,
  $\left(\begin{array}{cccccccc} 1&0&0&0&0&0&0&0 \\ 1&1&0&0&-1&0&1&0 \\ 1&1&0&0&0&-1&1&0 \\
        2&2&0&0&-1&-1&1&1 \\ 2&2&0&-1&0&-1&1&1 \\ 1&2&0&-1&0&0&0&1 \\
        0&1&1&-1&0&0&0&1 \\ 0&1&1&-1&0&0&0&0 \end{array}\right)$,
  $\left(\begin{array}{cccccccc} 1&0&0&0&0&0&0&0 \\ 1&1&0&0&-1&0&1&0 \\ 1&1&0&0&0&-1&1&0 \\
        2&2&0&0&-1&-1&2&-1 \\ 2&2&0&-1&0&-1&2&-1 \\ 1&2&0&-1&0&0&1&-1 \\
        0&1&1&-1&0&0&1&-1 \\ 0&0&0&0&0&0&1&-1 \end{array}\right)$,
        \\
  $\left(\begin{array}{cccccccc} 1&1&-1&-1&0&2&-1&0 \\ 0&1&0&-2&0&3&-1&0 \\ 1&2&-1&-2&-1&4&-1&0 \\
        1&2&-1&-3&-1&6&-2&1 \\ 1&1&-1&-2&-1&5&-2&1 \\
        1&1&-1&-2&0&3&-1&1 \\ 0&1&-1&-1&0&2&-1&1 \\ 0&1&0&-1&0&1&-1&1 \end{array}\right)$,
  $\left(\begin{array}{cccccccc} 1&1&-1&-1&0&2&-1&0 \\ 0&1&0&-2&0&3&-1&0 \\ 1&2&-1&-2&-1&4&-1&0 \\
        1&2&-1&-3&-1&6&-1&-1 \\ 1&1&-1&-2&-1&5&-1&-1 \\
        1&1&-1&-2&0&3&0&-1 \\ 0&1&-1&-1&0&2&0&-1 \\ 0&0&-1&0&0&1&0&0 \end{array}\right)$,
  $\left(\begin{array}{cccccccc} 1&1&-1&-1&0&2&-1&0 \\ 0&1&0&-2&0&3&-1&0 \\ 1&2&-1&-2&-1&4&-1&0 \\
        1&3&0&-4&-1&6&-2&0 \\ 1&2&0&-3&-1&5&-2&0 \\ 1&2&0&-3&0&3&-1&0 \\
        0&2&0&-2&0&2&-1&0 \\ 0&1&0&-1&0&1&0&-1 \end{array}\right)$,
        \\
  $\left(\begin{array}{cccccccc} 1&1&-1&-1&0&2&-1&0 \\ 0&2&-1&-2&1&2&-1&0 \\ 0&2&-1&-2&0&4&-2&0 \\
        0&3&-2&-3&1&5&-3&1 \\ 0&2&-2&-2&1&4&-3&1 \\ 0&1&-1&-2&1&3&-2&1 \\
        0&1&-1&-1&0&2&-1&1 \\ 0&1&0&-1&0&1&-1&1 \end{array}\right)$,
  $\left(\begin{array}{cccccccc} 1&1&-1&-1&0&2&-1&0 \\ 0&2&-1&-2&1&2&-1&0 \\ 0&2&-1&-2&0&4&-2&0 \\
        0&3&-2&-3&1&5&-2&-1 \\ 0&2&-2&-2&1&4&-2&-1 \\
        0&1&-1&-2&1&3&-1&-1 \\ 0&1&-1&-1&0&2&0&-1 \\ 0&0&-1&0&0&1&0&0 \end{array}\right)$,
  $\left(\begin{array}{cccccccc} 1&1&-1&-1&0&2&-1&0 \\ 0&2&-1&-2&1&2&-1&0 \\
        0&2&-1&-2&0&4&-2&0 \\ 0&4&-1&-4&1&5&-3&0 \\ 0&3&-1&-3&1&4&-3&0 \\
        0&2&0&-3&1&3&-2&0 \\ 0&2&0&-2&0&2&-1&0 \\ 0&1&0&-1&0&1&0&-1
        \end{array}\right)$,  \\
  $\left(\begin{array}{cccccccc} 1&1&-1&-1&0&2&-1&0 \\ 1&2&-1&-2&0&2&0&0 \\ 1&3&-2&-2&0&3&-1&0 \\
        2&4&-3&-3&0&4&-1&1 \\ 2&3&-3&-2&0&3&-1&1 \\ 1&2&-2&-2&1&2&-1&1 \\
        0&1&-1&-1&0&2&-1&1 \\ 0&1&0&-1&0&1&-1&1 \end{array}\right)$,
  $\left(\begin{array}{cccccccc} 1&1&-1&-1&0&2&-1&0 \\ 1&2&-1&-2&0&2&0&0 \\ 1&3&-2&-2&0&3&-1&0 \\
        2&4&-3&-3&0&4&0&-1 \\ 2&3&-3&-2&0&3&0&-1 \\ 1&2&-2&-2&1&2&0&-1 \\
        0&1&-1&-1&0&2&0&-1 \\ 0&0&-1&0&0&1&0&0 \end{array}\right)$,
  $\left(\begin{array}{cccccccc} 1&1&-1&-1&0&2&-1&0 \\ 1&2&-1&-2&0&2&0&0 \\ 1&3&-2&-2&0&3&-1&0 \\
        2&5&-2&-4&0&4&-1&0 \\ 2&4&-2&-3&0&3&-1&0 \\ 1&3&-1&-3&1&2&-1&0 \\
        0&2&0&-2&0&2&-1&0 \\ 0&1&0&-1&0&1&0&-1 \end{array}\right)$.
%%\end {document}
The character table of $G^{s_{95}}$:\\
${}_{
 \!\!\!\!_{
\!\!\!\!_{ \!\!\!\!  _{  \small
 \noindent
\begin{tabular}{c|ccccccccccccccccccccccc}
  & & & & & & & & & &10 & & & & & & & & & &20 & & &\\\hline
$\chi_{95}^{(1)}$&1&1&1&1&1&1&1&1&1&1&1&1&1&1&1 &1&1&1&1&1&1&1&1

\\  $\chi_{95}^{(2)}$&1&-1&-1&-1&-1&-1&-1&-1&-1&1&1&1&1&1&1&1
&1&1&1&1&1&1&1

\\  $\chi_{95}^{(3)}$&1&1&1&1&1&1&1&1&1&-1&-1&-1&-1&-1&-1&-1
&-1&-1&-1&-1&-1&-1&-1

\\  $\chi_{95}^{(4)}$&1&-1&-1&-1&-1&-1&-1&-1&-1&-1&-1&-1&-1&-1&-1&-1
&-1&-1&-1&-1&-1&-1&-1

\\  $\chi_{95}^{(5)}$&1&A&A&-1&-1&-1&/A&/A&/A&-A&-A&-A&1&1&1&-/A
&-/A&-/A&-A&-A&-A&1&1

\\  $\chi_{95}^{(6)}$&1&/A&/A&-1&-1&-1&A&A&A&-/A&-/A&-/A&1&1&1&-A
&-A&-A&-/A&-/A&-/A&1&1

\\  $\chi_{95}^{(7)}$&1&-/A&-/A&1&1&1&-A&-A&-A&/A&/A&/A&-1&-1&-1&A
&A&A&/A&/A&/A&-1&-1

\\
\end {tabular}
 }}}}$

${}_{
 \!\!\!\!\!\!\!_{
\!\!\!\!\!\!\!_{ \!\!\!\!\!\!\!  _{  \small
 \noindent
\begin{tabular}{c|ccccccccccccccccccccccc}
  & & & & & & & & & &10 & & & & & & & & & &20 & & &\\\hline

\\  $\chi_{95}^{({46})}$&1&-C&-C&-B&-B&-B&/C&/C&/C&-/A&-/A&-/A&1&1&1&-A
&-A&-A&/A&/A&/A&-1&-1

\\  $\chi_{95}^{({47})}$&1&-/C&-/C&B&B&B&C&C&C&-A&-A&-A&1&1&1&-/A
&-/A&-/A&A&A&A&-1&-1

\\  $\chi_{95}^{({48})}$&1&C&C&B&B&B&-/C&-/C&-/C&-/A&-/A&-/A&1&1&1&-A
&-A&-A&/A&/A&/A&-1&-1

\\  $\chi_{95}^{({49})}$&1&C&-/C&B&C&-/C&B&C&-/C&A&/A&-1&A&/A&-1&A
&/A&-1&-A&-/A&1&-A&-/A

\\  $\chi_{95}^{({50})}$&1&-/C&C&B&-/C&C&B&-/C&C&/A&A&-1&/A&A&-1&/A
&A&-1&-/A&-A&1&-/A&-A

\\  $\chi_{95}^{({51})}$&1&/C&-C&-B&/C&-C&-B&/C&-C&/A&A&-1&/A&A&-1&/A
&A&-1&-/A&-A&1&-/A&-A

\\
\end{tabular}
 }}}}$

\noindent \begin{tabular}{c|ccccccccccccccccccccccc}
  & & & & & & & 30& & & & & & & & & & 40& & & & & &\\\hline
$\chi_{95}^{(1)}$&1&1&1&1&1&1&1&1&1&1&1&1&1&1&1 &1&1&1&1&1&1&1&1

\\  $\chi_{95}^{(2)}$&1&1&1&1&-1&-1&-1&-1&-1&-1&-1&-1&-1&-1&-1
&-1&-1&-1&-1&-1&-1&-1&1

\\  $\chi_{95}^{(3)}$&-1&-1&-1&-1&1&1&1&1&1&1&1&1&1&-1&-1
&-1&-1&-1&-1&-1&-1&-1&1

\\  $\chi_{95}^{(4)}$&-1&-1&-1&-1&-1&-1&-1&-1&-1&-1&-1&-1&-1&1&1
&1&1&1&1&1&1&1&1

\\  $\chi_{95}^{(5)}$&1&-/A&-/A&-/A&A&A&A&-1&-1&-1&/A&/A&/A&/A&/A
&/A&-1&-1&-1&A&A&A&-/A

\\  $\chi_{95}^{(6)}$&1&-A&-A&-A&/A&/A&/A&-1&-1&-1&A&A&A&A&A
&A&-1&-1&-1&/A&/A&/A&-A

\\  $\chi_{95}^{(7)}$&-1&A&A&A&-/A&-/A&-/A&1&1&1&-A&-A&-A&A&A
&A&-1&-1&-1&/A&/A&/A&-A

\\  $\chi_{95}^{(8)}$&-1&/A&/A&/A&-A&-A&-A&1&1&1&-/A&-/A&-/A&/A&/A
&/A&-1&-1&-1&A&A&A&-/A

\\  $\chi_{95}^{(9)}$&-1&A&A&A&/A&/A&/A&-1&-1&-1&A&A&A&-A&-A
&-A&1&1&1&-/A&-/A&-/A&-A

\\  $\chi_{95}^{({10})}$&-1&/A&/A&/A&A&A&A&-1&-1&-1&/A&/A&/A&-/A&-/A
&-/A&1&1&1&-A&-A&-A&-/A

\\  $\chi_{95}^{({11})}$&1&-/A&-/A&-/A&-A&-A&-A&1&1&1&-/A&-/A&-/A&-/A&-/A
&-/A&1&1&1&-A&-A&-A&-/A

\\  $\chi_{95}^{({12})}$&1&-A&-A&-A&-/A&-/A&-/A&1&1&1&-A&-A&-A&-A&-A
&-A&1&1&1&-/A&-/A&-/A&-A

\\  $\chi_{95}^{({13})}$&1&-A&-/A&1&-1&/A&A&-1&/A&A&-1&/A&A&A&/A
&-1&A&/A&-1&A&/A&-1&1

\\  $\chi_{95}^{({14})}$&1&-/A&-A&1&-1&A&/A&-1&A&/A&-1&A&/A&/A&A
&-1&/A&A&-1&/A&A&-1&1

\\  $\chi_{95}^{({15})}$&1&1&-A&-/A&A&-1&/A&-1&/A&A&/A&A&-1&-1&A
&/A&A&/A&-1&/A&-1&A&-/A
\end{tabular}

${}_{
 \!\!\!\!\!\!\!_{
\!\!\!\!\!\!\!_{ \!\!\!\!\!\!\!  _{  \small
 \noindent
\begin{tabular}{c|ccccccccccccccccccccccc}
  & & & & & & & 30& & & & & & & & & & 40& & & & & &\\\hline
\\  $\chi_{95}^{({16})}$&1&1&-/A&-A&/A&-1&A&-1&A&/A&A&/A&-1&-1&/A
&A&/A&A&-1&A&-1&/A&-A

\\  $\chi_{95}^{({17})}$&1&-/A&1&-A&/A&A&-1&-1&/A&A&A&-1&/A&/A&-1
&A&A&/A&-1&-1&A&/A&-A

\\  $\chi_{95}^{({18})}$&1&-A&1&-/A&A&/A&-1&-1&A&/A&/A&-1&A&A&-1
&/A&/A&A&-1&-1&/A&A&-/A

\\  $\chi_{95}^{({19})}$&-1&A&/A&-1&1&-/A&-A&1&-/A&-A&1&-/A&-A&A&/A
&-1&A&/A&-1&A&/A&-1&1

\\  $\chi_{95}^{({20})}$&-1&/A&A&-1&1&-A&-/A&1&-A&-/A&1&-A&-/A&/A&A
&-1&/A&A&-1&/A&A&-1&1

\\  $\chi_{95}^{({21})}$&-1&/A&-1&A&-/A&-A&1&1&-/A&-A&-A&1&-/A&/A&-1
&A&A&/A&-1&-1&A&/A&-A

\\
\end{tabular}
 }}}}$

\noindent \begin{tabular}{c|ccccccccccccccccccccccc}
  & & & & & & & 30& & & & & & & & & & 40& & & & & &\\\hline
\\  $\chi_{95}^{({61})}$&-1&A&/A&-1&B&C&-/C&B&C&-/C&B&C&-/C&-/C&C
&B&-/C&C&B&-/C&C&B&-1

\\  $\chi_{95}^{({62})}$&-1&/A&A&-1&B&-/C&C&B&-/C&C&B&-/C&C&C&-/C
&B&C&-/C&B&C&-/C&B&-1

\\  $\chi_{95}^{({63})}$&-1&/A&A&-1&-B&/C&-C&-B&/C&-C&-B&/C&-C&-C&/C
&-B&-C&/C&-B&-C&/C&-B&-1

\\  $\chi_{95}^{({64})}$&-1&A&/A&-1&-B&-C&/C&-B&-C&/C&-B&-C&/C&/C&-C
&-B&/C&-C&-B&/C&-C&-B&-1

\\  $\chi_{95}^{({65})}$&-1&-1&A&/A&-/C&B&C&B&C&-/C&C&-/C&B&B&-/C
&C&-/C&C&B&C&B&-/C&/A

\\  $\chi_{95}^{({66})}$&-1&-1&/A&A&C&B&-/C&B&-/C&C&-/C&C&B&B&C
&-/C&C&-/C&B&-/C&B&C&A

\\  $\chi_{95}^{({67})}$&-1&-1&/A&A&-C&-B&/C&-B&/C&-C&/C&-C&-B&-B&-C
&/C&-C&/C&-B&/C&-B&-C&A

\\  $\chi_{95}^{({68})}$&-1&-1&A&/A&/C&-B&-C&-B&-C&/C&-C&/C&-B&-B&/C
&-C&/C&-C&-B&-C&-B&/C&/A

\\  $\chi_{95}^{({69})}$&-1&/A&-1&A&C&-/C&B&B&C&-/C&-/C&B&C&C&B
&-/C&-/C&C&B&B&-/C&C&A

\\  $\chi_{95}^{({70})}$&-1&A&-1&/A&-/C&C&B&B&-/C&C&C&B&-/C&-/C&B
&C&C&-/C&B&B&C&-/C&/A

\\  $\chi_{95}^{({71})}$&-1&A&-1&/A&/C&-C&-B&-B&/C&-C&-C&-B&/C&/C&-B
&-C&-C&/C&-B&-B&-C&/C&/A

\\  $\chi_{95}^{({72})}$&-1&/A&-1&A&-C&/C&-B&-B&-C&/C&/C&-B&-C&-C&-B
&/C&/C&-C&-B&-B&/C&-C&A\\
\end{tabular}

\noindent \begin{tabular}{c|cccccccccccccccccccccc}
  & & & &50 & & & & & & & & & &60 & & & & & & & &\\\hline
$\chi_{95}^{(1)}$&1&1&1&1&1&1&1&1&1&1&1&1&1&1&1 &1&1&1&1&1&1&1

\\  $\chi_{95}^{(2)}$&1&1&1&1&1&1&1&1&1&1&1&-1&1&1&1
&1&1&-1&-1&-1&-1&-1

\\  $\chi_{95}^{(3)}$&1&1&1&1&1&1&1&1&1&1&1&1&1&1&1
&1&1&-1&-1&-1&-1&-1

\\  $\chi_{95}^{(4)}$&1&1&1&1&1&1&1&1&1&1&1&-1&1&1&1
&1&1&1&1&1&1&1

\\  $\chi_{95}^{(5)}$&-/A&-/A&1&1&1&-A&-A&-A&-/A&-/A&-/A&A&1&1&-A
&-A&-A&/A&/A&/A&-1&-1

\\  $\chi_{95}^{(6)}$&-A&-A&1&1&1&-/A&-/A&-/A&-A&-A&-A&/A&1&1&-/A
&-/A&-/A&A&A&A&-1&-1

\\  $\chi_{95}^{(7)}$&-A&-A&1&1&1&-/A&-/A&-/A&-A&-A&-A&-/A&1&1&-/A
&-/A&-/A&A&A&A&-1&-1

\\
\end{tabular}

\noindent \begin{tabular}{c|cccccccccccccccccccccc}
  & & & &50 & & & & & & & & & &60 & & & & & & & &\\\hline

\\  $\chi_{95}^{({31})}$&-A&-/A&1&-A&-/A&1&-A&-/A&1&-A&-/A&1&-A&-/A&1
&-A&-/A&-/A&-A&1&-/A&-A

\\  $\chi_{95}^{({32})}$&-/A&-A&1&-/A&-A&1&-/A&-A&1&-/A&-A&1&-/A&-A&1
&-/A&-A&-A&-/A&1&-A&-/A

\\  $\chi_{95}^{({33})}$&1&-A&1&-A&-/A&-A&-/A&1&-/A&1&-A&-A&-A&-/A&-A
&-/A&1&-A&1&-/A&-/A&-A

\\  $\chi_{95}^{({34})}$&1&-/A&1&-/A&-A&-/A&-A&1&-A&1&-/A&-/A&-/A&-A&-/A
&-A&1&-/A&1&-A&-A&-/A

\\  $\chi_{95}^{({35})}$&-/A&1&1&-A&-/A&-/A&1&-A&-A&-/A&1&-/A&-A&-/A&-/A
&1&-A&1&-/A&-A&-/A&-A

\\  $\chi_{95}^{({36})}$&-A&1&1&-/A&-A&-A&1&-/A&-/A&-A&1&-A&-/A&-A&-A
&1&-/A&1&-A&-/A&-A&-/A

\\
\end{tabular}

\noindent \begin{tabular}{c|cccc}
  & &70 & &\\\hline
$\chi_{95}^{(1)}$&1&1&1&1

\\  $\chi_{95}^{(2)}$&-1&-1&-1&-1

\\  $\chi_{95}^{(3)}$&-1&-1&-1&-1

\\  $\chi_{95}^{(4)}$&1&1&1&1

\\  $\chi_{95}^{(5)}$&-1&A&A&A

\\  $\chi_{95}^{(6)}$&-1&/A&/A&/A

\\  $\chi_{95}^{(7)}$&-1&/A&/A&/A

\\  $\chi_{95}^{(8)}$&-1&A&A&A

\\  $\chi_{95}^{(9)}$&1&-/A&-/A&-/A

\\  $\chi_{95}^{({10})}$&1&-A&-A&-A

\\  $\chi_{95}^{({11})}$&1&-A&-A&-A

\\  $\chi_{95}^{({12})}$&1&-/A&-/A&-/A

\\  $\chi_{95}^{({13})}$&-1&A&/A&-1

\\  $\chi_{95}^{({14})}$&-1&/A&A&-1

\\  $\chi_{95}^{({15})}$&-1&/A&-1&A

\\  $\chi_{95}^{({16})}$&-1&A&-1&/A

\\  $\chi_{95}^{({17})}$&-1&-1&A&/A

\\  $\chi_{95}^{({18})}$&-1&-1&/A&A

\\  $\chi_{95}^{({19})}$&-1&A&/A&-1

\\  $\chi_{95}^{({20})}$&-1&/A&A&-1

\\  $\chi_{95}^{({21})}$&-1&-1&A&/A

\\  $\chi_{95}^{({22})}$&-1&-1&/A&A

\\  $\chi_{95}^{({23})}$&-1&/A&-1&A

\\  $\chi_{95}^{({24})}$&-1&A&-1&/A

\\  $\chi_{95}^{({25})}$&1&-/A&-A&1

\\  $\chi_{95}^{({26})}$&1&-A&-/A&1

\\  $\chi_{95}^{({27})}$&1&-A&1&-/A

\\  $\chi_{95}^{({28})}$&1&-/A&1&-A

\\  $\chi_{95}^{({29})}$&1&1&-/A&-A

\\  $\chi_{95}^{({30})}$&1&1&-A&-/A

\\  $\chi_{95}^{({31})}$&1&-/A&-A&1

\\  $\chi_{95}^{({32})}$&1&-A&-/A&1

\\  $\chi_{95}^{({33})}$&1&1&-/A&-A

\\  $\chi_{95}^{({34})}$&1&1&-A&-/A

\\  $\chi_{95}^{({35})}$&1&-A&1&-/A

\\  $\chi_{95}^{({36})}$&1&-/A&1&-A

\\  $\chi_{95}^{({37})}$&-B&-B&-B&-B

\\  $\chi_{95}^{({38})}$&B&B&B&B

\\  $\chi_{95}^{({39})}$&-B&-C&-C&-C

\\  $\chi_{95}^{({40})}$&-B&/C&/C&/C

\\  $\chi_{95}^{({41})}$&B&C&C&C

\\  $\chi_{95}^{({42})}$&B&-/C&-/C&-/C

\\  $\chi_{95}^{({43})}$&-B&-B&-B&-B

\\  $\chi_{95}^{({44})}$&B&B&B&B

\\  $\chi_{95}^{({45})}$&-B&/C&/C&/C
\end{tabular}

\noindent \begin{tabular}{c|cccc}
  & &70 & &\\\hline

\\  $\chi_{95}^{({46})}$&-B&-C&-C&-C

\\  $\chi_{95}^{({47})}$&B&-/C&-/C&-/C

\\  $\chi_{95}^{({48})}$&B&C&C&C

\\  $\chi_{95}^{({49})}$&-B&/C&-C&-B

\\  $\chi_{95}^{({50})}$&-B&-C&/C&-B

\\  $\chi_{95}^{({51})}$&B&C&-/C&B

\\  $\chi_{95}^{({52})}$&B&-/C&C&B

\\  $\chi_{95}^{({53})}$&-B&-B&/C&-C

\\  $\chi_{95}^{({54})}$&-B&-B&-C&/C

\\  $\chi_{95}^{({55})}$&B&B&C&-/C

\\  $\chi_{95}^{({56})}$&B&B&-/C&C

\\  $\chi_{95}^{({57})}$&-B&-C&-B&/C

\\  $\chi_{95}^{({58})}$&-B&/C&-B&-C

\\  $\chi_{95}^{({59})}$&B&-/C&B&C

\\  $\chi_{95}^{({60})}$&B&C&B&-/C

\\  $\chi_{95}^{({61})}$&-B&/C&-C&-B

\\  $\chi_{95}^{({62})}$&-B&-C&/C&-B

\\  $\chi_{95}^{({63})}$&B&C&-/C&B

\\  $\chi_{95}^{({64})}$&B&-/C&C&B

\\  $\chi_{95}^{({65})}$&-B&-C&-B&/C

\\  $\chi_{95}^{({66})}$&-B&/C&-B&-C

\\  $\chi_{95}^{({67})}$&B&-/C&B&C

\\  $\chi_{95}^{({68})}$&B&C&B&-/C

\\  $\chi_{95}^{({69})}$&-B&-B&/C&-C

\\  $\chi_{95}^{({70})}$&-B&-B&-C&/C

\\  $\chi_{95}^{({71})}$&B&B&C&-/C

\\  $\chi_{95}^{({72})}$&B&B&-/C&C\\
\end{tabular}

\noindent where A = -E(3)
  = (1-ER(-3))/2 = -b3,
B = E(4)
  = ER(-1) = i,
C = E(12)$^{11}$.

The generators of $G^{s_{96}}$ are:\\
$\left(\begin{array}{cccccccc} 1&-1&0&-1& 2&-1&0&0 \\
1&-1&0&-1&2&-2&1&0 \\ 1&-1&0&-2&4&-3&
1&0 \\
        2&-1&0&-3&5&-4&2&-1 \\ 2&-1&0&-2&4&-4&2&-1 \\
        1&-1&0&-1&3&-3&1&-1 \\ 0&-1&1&-1&2&-2&1&-1 \\
        0&-1&0&0&1&-1&1&-1 \end{array}\right)$,
  $\left(\begin{array}{cccccccc} 1&0&0&-2&2&0&0&0 \\ 0&1&0&-3&3&0&0&0 \\ 0&0&1&-4&4&0&0&0 \\
        0&0&0&-5&6&0&0&0 \\ 0&0&0&-4&5&0&0&0 \\ 0&0&0&-3&3&1&0&0 \\
        0&0&0&-2&2&0&1&0 \\ 0&0&0&-1&1&0&0&1 \end{array}\right)$,
  $\left(\begin{array}{cccccccc} 1&1&-1&0&-1&2&-1&0 \\ 0&2&-1&0&-1&2&-1&0 \\ 0&2&-1&0&-2&4&-2&0 \\
        0&3&-2&0&-2&5&-3&1 \\ 0&2&-2&1&-2&4&-3&1 \\ 0&1&-1&1&-2&3&-2&1 \\
        0&1&-1&1&-2&2&-1&1 \\ 0&1&0&0&-1&1&-1&1 \end{array}\right)$,
  $\left(\begin{array}{cccccccc} -1&1&1&-1&0&0&1&0 \\ 0&1&1&-1&-1&1&1&0 \\ 0&2&1&-2&0&0&2&0 \\
        0&2&1&-2&-1&1&3&0 \\ 0&1&1&-1&-1&0&3&0 \\ 0&1&0&0&-1&0&2&0 \\
        0&0&0&0&0&0&1&0 \\ 0&0&0&0&0&0&0&1 \end{array}\right)$,
  $\left(\begin{array}{cccccccc} -1&0&0&2&-2&0&0&0 \\ -1&0&1&2&-2&-1&0&0 \\ -2&0&1&3&-3&-1&1&-1 \\
        -2&0&2&4&-5&-1&1&-1 \\ -2&0&2&3&-4&-1&1&-1 \\
        -1&-1&1&3&-3&-1&1&-1 \\ -1&-1&1&2&-2&-1&1&0 \\
        0&0&0&1&-1&-1&1&0 \end{array}\right)$.

The representatives of conjugacy classes of   $G^{s_{96}}$ are:\\
$\left(\begin{array}{cccccccc} 1&0&0&0&0&0&0&0 \\ 0&1&0&0&
0&0&0&0 \\ 0&0&1&0&0&0&0&0 \\
        0&0&0&1&0&0&0&0 \\ 0&0&0&0&1&0&0&0 \\ 0&0&0&0&0&1&0&0 \\
        0&0&0&0&0&0&1&0 \\ 0&0&0&0&0&0&0&1 \end{array}\right)$,
  $\left(\begin{array}{cccccccc} -1&-1&1&0&1&-2&1&0 \\ -1&-2&0&1&1&-2&1&0 \\ -2&-2&1&1&1&-3&1&1 \\
        -2&-4&1&2&1&-4&2&1 \\ -2&-3&1&1&1&-3&2&1 \\ -1&-2&1&0&1&-2&1&1 \\
        -1&-1&1&0&0&-1&1&0 \\ 0&0&0&0&0&-1&1&0 \end{array}\right)$,
  $\left(\begin{array}{cccccccc} -1&-1&1&0&1&-2&1&0 \\ -1&-2&0&1&1&-2&1&0 \\ -2&-2&1&1&1&-3&1&1 \\
        -2&-3&0&2&2&-5&2&1 \\ -2&-2&0&1&2&-4&2&1 \\ -1&-1&0&0&2&-3&1&1 \\
        -1&0&0&0&1&-2&1&0 \\ -1&0&0&0&1&-1&0&0 \end{array}\right)$,
        \\
  $\left(\begin{array}{cccccccc} -1&-1&1&0&1&-2&1&0 \\ -1&-2&1&0&2&-2&0&0 \\ -1&-3&2&0&2&-3&1&0 \\
        -2&-5&2&1&3&-4&1&0 \\ -2&-4&2&0&3&-3&1&0 \\ -1&-3&1&0&2&-2&1&0 \\
        0&-2&0&0&2&-2&1&0 \\ 0&-1&0&0&1&-1&0&1 \end{array}\right)$,
  $\left(\begin{array}{cccccccc} -1&-1&1&0&1&-2&1&0 \\ -1&-2&1&0&2&-2&0&0 \\ -1&-3&2&0&2&-3&1&0 \\
        -2&-4&3&0&3&-4&0&1 \\ -2&-3&3&-1&3&-3&0&1 \\
        -1&-2&2&-1&2&-2&0&1 \\ 0&-1&1&-1&2&-2&0&1 \\ 0&0&1&-1&1&-1&0&0 \end{array}\right)$,
 $\left(\begin{array}{cccccccc} -1&-1&1&0&1&-2&1&0 \\ -1&-2&1&0&2&-2&0&0 \\
        -1&-3&2&0&2&-3&1&0 \\ -2&-4&3&0&3&-4&1&-1 \\
        -2&-3&3&-1&3&-3&1&-1 \\ -1&-2&2&-1&2&-2&1&-1 \\
        0&-1&1&-1&2&-2&1&-1 \\ 0&-1&0&0&1&-1&1&-1
        \end{array}\right)$,  \\
  $\left(\begin{array}{cccccccc} -1&-1&1&0&1&-2&1&0 \\ 0&-2&1&0&1&-2&1&0 \\ 0&-2&1&0&2&-4&2&0 \\
        0&-4&1&1&2&-5&3&0 \\ 0&-3&1&0&2&-4&3&0 \\ 0&-2&0&0&2&-3&2&0 \\
        0&-2&0&0&2&-2&1&0 \\ 0&-1&0&0&1&-1&0&1 \end{array}\right)$,
  $\left(\begin{array}{cccccccc} -1&-1&1&0&1&-2&1&0 \\ 0&-2&1&0&1&-2&1&0 \\ 0&-2&1&0&2&-4&2&0 \\
        0&-3&2&0&2&-5&2&1 \\ 0&-2&2&-1&2&-4&2&1 \\ 0&-1&1&-1&2&-3&1&1 \\
        0&-1&1&-1&2&-2&0&1 \\ 0&0&1&-1&1&-1&0&0 \end{array}\right)$,
  $\left(\begin{array}{cccccccc} -1&-1&1&0&1&-2&1&0 \\ 0&-2&1&0&1&-2&1&0 \\ 0&-2&1&0&2&-4&2&0 \\
        0&-3&2&0&2&-5&3&-1 \\ 0&-2&2&-1&2&-4&3&-1 \\
        0&-1&1&-1&2&-3&2&-1 \\ 0&-1&1&-1&2&-2&1&-1 \\
        0&-1&0&0&1&-1&1&-1 \end{array}\right)$,  \\
  $\left(\begin{array}{cccccccc} -1&-1&1&0&1&-2&1&0 \\ 0&-1&0&0&2&-3&1&0 \\ -1&-2&1&0&3&-4&1&0 \\
        -1&-3&0&1&4&-6&2&0 \\ -1&-2&0&0&4&-5&2&0 \\ -1&-2&0&0&3&-3&1&0 \\
        0&-2&0&0&2&-2&1&0 \\ 0&-1&0&0&1&-1&0&1 \end{array}\right)$,
  $\left(\begin{array}{cccccccc} -1&-1&1&0&1&-2&1&0 \\ 0&-1&0&0&2&-3&1&0 \\ -1&-2&1&0&3&-4&1&0 \\
        -1&-2&1&0&4&-6&1&1 \\ -1&-1&1&-1&4&-5&1&1 \\
        -1&-1&1&-1&3&-3&0&1 \\ 0&-1&1&-1&2&-2&0&1 \\ 0&0&1&-1&1&-1&0&0 \end{array}\right)$,
  $\left(\begin{array}{cccccccc} -1&-1&1&0&1&-2&1&0 \\ 0&-1&0&0&2&-3&1&0 \\ -1&-2&1&0&3&-4&1&0 \\
        -1&-2&1&0&4&-6&2&-1 \\ -1&-1&1&-1&4&-5&2&-1 \\
        -1&-1&1&-1&3&-3&1&-1 \\ 0&-1&1&-1&2&-2&1&-1 \\
        0&-1&0&0&1&-1&1&-1 \end{array}\right)$,  \\
  $\left(\begin{array}{cccccccc} -1&-1&1&1&0&-2&1&0 \\ -1&-2&0&3&-1&-2&1&0 \\ -2&-2&1&3&-1&-3&1&1 \\
        -3&-4&1&5&-1&-4&1&1 \\ -3&-3&1&4&-1&-3&1&1 \\
        -2&-2&1&3&-1&-2&0&1 \\ -2&-1&1&2&-1&-1&0&0 \\
        -1&-1&1&1&-1&0&0&0 \end{array}\right)$,
  $\left(\begin{array}{cccccccc} -1&-1&1&1&0&-2&1&0 \\ -1&-2&0&3&-1&-2&1&0 \\ -2&-2&1&3&-1&-3&1&1 \\
        -2&-4&1&5&-2&-4&2&1 \\ -2&-3&1&4&-2&-3&2&1 \\
        -1&-2&1&3&-2&-2&1&1 \\ -1&-1&1&2&-2&-1&1&0 \\
        0&0&0&1&-1&-1&1&0 \end{array}\right)$,
  $\left(\begin{array}{cccccccc} -1&-1&1&1&0&-2&1&0 \\ -1&-2&0&3&-1&-2&1&0 \\ -2&-2&1&3&-1&-3&1&1 \\
        -2&-3&0&5&-1&-5&2&1 \\ -2&-2&0&4&-1&-4&2&1 \\
        -1&-1&0&3&-1&-3&1&1 \\ -1&0&0&2&-1&-2&1&0 \\ -1&0&0&1&0&-1&0&0
        \end{array}\right)$,  \\
  $\left(\begin{array}{cccccccc} -1&-1&1&1&0&-2&1&0 \\ -1&-2&1&2&0&-2&0&0 \\
        -1&-3&2&2&0&-3&1&0 \\ -2&-5&2&4&0&-4&1&0 \\ -2&-4&2&3&0&-3&1&0 \\
        -1&-3&1&3&-1&-2&1&0 \\ 0&-2&0&2&0&-2&1&0 \\ 0&-1&0&1&0&-1&0&1 \end{array}\right)$,
  $\left(\begin{array}{cccccccc} -1&-1&1&1&0&-2&1&0 \\ -1&-2&1&2&0&-2&0&0 \\
        -1&-3&2&2&0&-3&1&0 \\ -2&-4&3&3&0&-4&0&1 \\ -2&-3&3&2&0&-3&0&1 \\
        -1&-2&2&2&-1&-2&0&1 \\ 0&-1&1&1&0&-2&0&1 \\ 0&0&1&0&0&-1&0&0 \end{array}\right)$,
  $\left(\begin{array}{cccccccc} -1&-1&1&1&0&-2&1&0 \\ -1&-2&1&2&0&-2&0&0 \\
        -1&-3&2&2&0&-3&1&0 \\ -2&-4&3&3&0&-4&1&-1 \\
        -2&-3&3&2&0&-3&1&-1 \\ -1&-2&2&2&-1&-2&1&-1 \\
        0&-1&1&1&0&-2&1&-1 \\ 0&-1&0&1&0&-1&1&-1
        \end{array}\right)$,  \\
  $\left(\begin{array}{cccccccc} -1&-1&1&1&0&-2&1&0 \\ 0&-2&1&2&-1&-2&1&0 \\ 0&-2&1&2&0&-4&2&0 \\
        0&-4&1&4&-1&-5&3&0 \\ 0&-3&1&3&-1&-4&3&0 \\ 0&-2&0&3&-1&-3&2&0 \\
        0&-2&0&2&0&-2&1&0 \\ 0&-1&0&1&0&-1&0&1 \end{array}\right)$,
  $\left(\begin{array}{cccccccc} -1&-1&1&1&0&-2&1&0 \\ 0&-2&1&2&-1&-2&1&0 \\ 0&-2&1&2&0&-4&2&0 \\
        0&-3&2&3&-1&-5&2&1 \\ 0&-2&2&2&-1&-4&2&1 \\ 0&-1&1&2&-1&-3&1&1 \\
        0&-1&1&1&0&-2&0&1 \\ 0&0&1&0&0&-1&0&0 \end{array}\right)$,
  $\left(\begin{array}{cccccccc} -1&-1&1&1&0&-2&1&0 \\ 0&-2&1&2&-1&-2&1&0 \\ 0&-2&1&2&0&-4&2&0 \\
        0&-3&2&3&-1&-5&3&-1 \\ 0&-2&2&2&-1&-4&3&-1 \\
        0&-1&1&2&-1&-3&2&-1 \\ 0&-1&1&1&0&-2&1&-1 \\
        0&-1&0&1&0&-1&1&-1 \end{array}\right)$,  \\
  $\left(\begin{array}{cccccccc} -1&-1&1&1&0&-2&1&0 \\ 0&-1&0&2&0&-3&1&0 \\ -1&-2&1&2&1&-4&1&0 \\
        -1&-3&0&4&1&-6&2&0 \\ -1&-2&0&3&1&-5&2&0 \\ -1&-2&0&3&0&-3&1&0 \\
        0&-2&0&2&0&-2&1&0 \\ 0&-1&0&1&0&-1&0&1 \end{array}\right)$,
  $\left(\begin{array}{cccccccc} -1&-1&1&1&0&-2&1&0 \\ 0&-1&0&2&0&-3&1&0 \\ -1&-2&1&2&1&-4&1&0 \\
        -1&-2&1&3&1&-6&1&1 \\ -1&-1&1&2&1&-5&1&1 \\ -1&-1&1&2&0&-3&0&1 \\
        0&-1&1&1&0&-2&0&1 \\ 0&0&1&0&0&-1&0&0 \end{array}\right)$,
  $\left(\begin{array}{cccccccc} -1&-1&1&1&0&-2&1&0 \\ 0&-1&0&2&0&-3&1&0 \\ -1&-2&1&2&1&-4&1&0 \\
        -1&-2&1&3&1&-6&2&-1 \\ -1&-1&1&2&1&-5&2&-1 \\
        -1&-1&1&2&0&-3&1&-1 \\ 0&-1&1&1&0&-2&1&-1 \\
        0&-1&0&1&0&-1&1&-1 \end{array}\right)$,  \\
  $\left(\begin{array}{cccccccc} -1&0&0&0&0&0&0&0 \\ -1&-1&0&0&1&-1&0&1 \\ -2&0&1&-1&1&-1&0&1 \\
        -3&-1&1&-1&2&-1&-1&2 \\ -3&-1&1&0&1&-1&-1&2 \\
        -2&-1&1&0&1&-1&-1&1 \\ -2&-1&1&0&1&-1&0&0 \\ -1&-1&1&0&0&0&0&0 \end{array}\right)$,
  $\left(\begin{array}{cccccccc} -1&0&0&0&0&0&0&0 \\ -1&-1&0&0&1&-1&0&1 \\ -2&0&1&-1&1&-1&0&1 \\
        -2&-1&1&-1&1&-1&0&2 \\ -2&-1&1&0&0&-1&0&2 \\
        -1&-1&1&0&0&-1&0&1 \\ -1&-1&1&0&0&-1&1&0 \\ 0&0&0&0&0&-1&1&0 \end{array}\right)$,
  $\left(\begin{array}{cccccccc} -1&0&0&0&0&0&0&0 \\ -1&-1&0&0&1&-1&0&1 \\ -2&0&1&-1&1&-1&0&1 \\
        -2&0&0&-1&2&-2&0&2 \\ -2&0&0&0&1&-2&0&2 \\ -1&0&0&0&1&-2&0&1 \\
        -1&0&0&0&1&-2&1&0 \\ -1&0&0&0&1&-1&0&0 \end{array}\right)$,
        \\
  $\left(\begin{array}{cccccccc} -1&0&0&0&0&0&0&0 \\ -1&-1&0&0&1&0&-1&0 \\ -1&-1&0&0&0&1&-1&0 \\
        -2&-2&0&0&1&1&-2&1 \\ -2&-2&0&1&0&1&-2&1 \\ -1&-2&0&1&0&0&-1&1 \\
        0&-1&-1&1&0&0&-1&1 \\ 0&0&0&0&0&0&-1&1 \end{array}\right)$,
  $\left(\begin{array}{cccccccc} -1&0&0&0&0&0&0&0 \\ -1&-1&0&0&1&0&-1&0 \\ -1&-1&0&0&0&1&-1&0 \\
        -2&-2&0&0&1&1&-1&-1 \\ -2&-2&0&1&0&1&-1&-1 \\
        -1&-2&0&1&0&0&0&-1 \\ 0&-1&-1&1&0&0&0&-1 \\ 0&-1&-1&1&0&0&0&0 \end{array}\right)$,
  $\left(\begin{array}{cccccccc} -1&0&0&0&0&0&0&0 \\ -1&-1&0&0&1&0&-1&0 \\ -1&-1&0&0&0&1&-1&0 \\
        -2&-1&1&-1&1&1&-2&0 \\ -2&-1&1&0&0&1&-2&0 \\
        -1&-1&1&0&0&0&-1&0 \\ 0&0&0&0&0&0&-1&0 \\ 0&0&0&0&0&0&0&-1
        \end{array}\right)$,  \\
  $\left(\begin{array}{cccccccc} -1&0&0&0&0&0&0&0 \\ 0&-1&0&0&0&0&0&0 \\ 0&0&-1&0&0&0&0&0 \\
        0&-1&-1&0&0&0&1&-1 \\ 0&-1&-1&1&-1&0&1&-1 \\
        0&-1&-1&1&0&-1&1&-1 \\ 0&-1&-1&1&0&0&0&-1 \\ 0&-1&-1&1&0&0&0&0 \end{array}\right)$,
  $\left(\begin{array}{cccccccc} -1&0&0&0&0&0&0&0 \\ 0&-1&0&0&0&0&0&0 \\ 0&0&-1&0&0&0&0&0 \\
        0&0&0&-1&0&0&0&0 \\ 0&0&0&0&-1&0&0&0 \\ 0&0&0&0&0&-1&0&0 \\
        0&0&0&0&0&0&-1&0 \\ 0&0&0&0&0&0&0&-1 \end{array}\right)$,
  $\left(\begin{array}{cccccccc} -1&0&0&0&0&0&0&0 \\ 0&0&-1&0&1&-1&0&0 \\ -1&0&-1&0&1&0&-1&0 \\
        -1&0&-2&0&2&-1&0&-1 \\ -1&0&-2&1&1&-1&0&-1 \\
        -1&-1&-1&1&1&-1&0&-1 \\ 0&-1&-1&1&0&0&0&-1 \\
        0&-1&-1&1&0&0&0&0 \end{array}\right)$,  \\
  $\left(\begin{array}{cccccccc} -1&0&0&2&-2&0&0&0 \\ -1&-1&0&3&-2&-1&0&1 \\ -2&0&1&3&-3&-1&0&1 \\
        -3&-1&1&5&-4&-1&-1&2 \\ -3&-1&1&4&-3&-1&-1&2 \\
        -2&-1&1&3&-2&-1&-1&1 \\ -2&-1&1&2&-1&-1&0&0 \\
        -1&-1&1&1&-1&0&0&0 \end{array}\right)$,
  $\left(\begin{array}{cccccccc} -1&0&0&2&-2&0&0&0 \\ -1&-1&0&3&-2&-1&0&1 \\ -2&0&1&3&-3&-1&0&1 \\
        -2&-1&1&5&-5&-1&0&2 \\ -2&-1&1&4&-4&-1&0&2 \\
        -1&-1&1&3&-3&-1&0&1 \\ -1&-1&1&2&-2&-1&1&0 \\
        0&0&0&1&-1&-1&1&0 \end{array}\right)$,
  $\left(\begin{array}{cccccccc} -1&0&0&2&-2&0&0&0 \\ -1&-1&0&3&-2&-1&0&1 \\ -2&0&1&3&-3&-1&0&1 \\
        -2&0&0&5&-4&-2&0&2 \\ -2&0&0&4&-3&-2&0&2 \\ -1&0&0&3&-2&-2&0&1 \\
        -1&0&0&2&-1&-2&1&0 \\ -1&0&0&1&0&-1&0&0 \end{array}\right)$,
         \\
  $\left(\begin{array}{cccccccc} -1&0&0&2&-2&0&0&0 \\ -1&-1&0&3&-2&0&-1&0 \\ -1&-1&0&4&-4&1&-1&0 \\
        -2&-2&0&6&-5&1&-2&1 \\ -2&-2&0&5&-4&1&-2&1 \\
        -1&-2&0&4&-3&0&-1&1 \\ 0&-1&-1&3&-2&0&-1&1 \\
        0&0&0&1&-1&0&-1&1 \end{array}\right)$,
  $\left(\begin{array}{cccccccc} -1&0&0&2&-2&0&0&0 \\ -1&-1&0&3&-2&0&-1&0 \\ -1&-1&0&4&-4&1&-1&0 \\
        -2&-2&0&6&-5&1&-1&-1 \\ -2&-2&0&5&-4&1&-1&-1 \\
        -1&-2&0&4&-3&0&0&-1 \\ 0&-1&-1&3&-2&0&0&-1 \\
        0&-1&-1&2&-1&0&0&0 \end{array}\right)$,
  $\left(\begin{array}{cccccccc} -1&0&0&2&-2&0&0&0 \\ -1&-1&0&3&-2&0&-1&0 \\ -1&-1&0&4&-4&1&-1&0 \\
        -2&-1&1&5&-5&1&-2&0 \\ -2&-1&1&4&-4&1&-2&0 \\
        -1&-1&1&3&-3&0&-1&0 \\ 0&0&0&2&-2&0&-1&0 \\ 0&0&0&1&-1&0&0&-1
        \end{array}\right)$,  \\
  $\left(\begin{array}{cccccccc} -1&0&0&2&-2&0&0&0 \\ 0&-1&0&3&-3&0&0&0 \\ 0&0&-1&4&-4&0&0&0 \\
        0&-1&-1&6&-6&0&1&-1 \\ 0&-1&-1&5&-5&0&1&-1 \\
        0&-1&-1&4&-3&-1&1&-1 \\ 0&-1&-1&3&-2&0&0&-1 \\
        0&-1&-1&2&-1&0&0&0 \end{array}\right)$,
  $\left(\begin{array}{cccccccc} -1&0&0&2&-2&0&0&0 \\ 0&-1&0&3&-3&0&0&0 \\ 0&0&-1&4&-4&0&0&0 \\
        0&0&0&5&-6&0&0&0 \\ 0&0&0&4&-5&0&0&0 \\ 0&0&0&3&-3&-1&0&0 \\
        0&0&0&2&-2&0&-1&0 \\ 0&0&0&1&-1&0&0&-1 \end{array}\right)$,
  $\left(\begin{array}{cccccccc} -1&0&0&2&-2&0&0&0 \\ 0&0&-1&3&-2&-1&0&0 \\ -1&0&-1&4&-3&0&-1&0 \\
        -1&0&-2&6&-4&-1&0&-1 \\ -1&0&-2&5&-3&-1&0&-1 \\
        -1&-1&-1&4&-2&-1&0&-1 \\ 0&-1&-1&3&-2&0&0&-1 \\
        0&-1&-1&2&-1&0&0&0 \end{array}\right)$,  \\
  $\left(\begin{array}{cccccccc} -1&1&0&-1&0&1&0&0 \\ -1&1&0&-2&1&1&0&1 \\ -2&2&1&-3&1&1&0&1 \\
        -3&2&1&-4&2&2&-1&2 \\ -3&2&1&-3&1&2&-1&2 \\ -2&2&1&-3&1&2&-1&1 \\
        -2&1&1&-2&1&1&0&0 \\ -1&0&1&-1&0&1&0&0 \end{array}\right)$,
  $\left(\begin{array}{cccccccc} -1&1&0&-1&0&1&0&0 \\ -1&1&0&-2&1&1&0&1 \\ -2&2&1&-3&1&1&0&1 \\
        -2&2&1&-4&1&2&0&2 \\ -2&2&1&-3&0&2&0&2 \\ -1&2&1&-3&0&2&0&1 \\
        -1&1&1&-2&0&1&1&0 \\ 0&1&0&-1&0&0&1&0 \end{array}\right)$,
  $\left(\begin{array}{cccccccc} -1&1&0&-1&0&1&0&0 \\ -1&1&0&-2&1&1&0&1 \\ -2&2&1&-3&1&1&0&1 \\
        -2&3&0&-4&2&1&0&2 \\ -2&3&0&-3&1&1&0&2 \\ -1&3&0&-3&1&1&0&1 \\
        -1&2&0&-2&1&0&1&0 \\ -1&1&0&-1&1&0&0&0 \end{array}\right)$,
        \\
  $\left(\begin{array}{cccccccc} -1&1&0&1&-2&1&0&0 \\ -1&1&0&1&-2&1&0&1 \\ -2&2&1&1&-3&1&0&1 \\
        -3&2&1&2&-4&2&-1&2 \\ -3&2&1&1&-3&2&-1&2 \\ -2&2&1&0&-2&2&-1&1 \\
        -2&1&1&0&-1&1&0&0 \\ -1&0&1&0&-1&1&0&0 \end{array}\right)$,
  $\left(\begin{array}{cccccccc} -1&1&0&1&-2&1&0&0 \\ -1&1&0&1&-2&1&0&1 \\ -2&2&1&1&-3&1&0&1 \\
        -2&2&1&2&-5&2&0&2 \\ -2&2&1&1&-4&2&0&2 \\ -1&2&1&0&-3&2&0&1 \\
        -1&1&1&0&-2&1&1&0 \\ 0&1&0&0&-1&0&1&0 \end{array}\right)$,
  $\left(\begin{array}{cccccccc} -1&1&0&1&-2&1&0&0 \\ -1&1&0&1&-2&1&0&1 \\ -2&2&1&1&-3&1&0&1 \\
        -2&3&0&2&-4&1&0&2 \\ -2&3&0&1&-3&1&0&2 \\ -1&3&0&0&-2&1&0&1 \\
        -1&2&0&0&-1&0&1&0 \\ -1&1&0&0&0&0&0&0 \end{array}\right)$,
        \\
  $\left(\begin{array}{cccccccc} -1&1&0&1&-2&1&0&0 \\ -1&1&0&1&-2&2&-1&0 \\ -1&1&0&2&-4&3&-1&0 \\
        -2&1&0&3&-5&4&-2&1 \\ -2&1&0&2&-4&4&-2&1 \\ -1&1&0&1&-3&3&-1&1 \\
        0&1&-1&1&-2&2&-1&1 \\ 0&1&0&0&-1&1&-1&1 \end{array}\right)$,
  $\left(\begin{array}{cccccccc} -1&1&0&1&-2&1&0&0 \\ -1&1&0&1&-2&2&-1&0 \\ -1&1&0&2&-4&3&-1&0 \\
        -2&1&0&3&-5&4&-1&-1 \\ -2&1&0&2&-4&4&-1&-1 \\
        -1&1&0&1&-3&3&0&-1 \\ 0&1&-1&1&-2&2&0&-1 \\ 0&0&-1&1&-1&1&0&0 \end{array}\right)$,
  $\left(\begin{array}{cccccccc} -1&1&0&1&-2&1&0&0 \\ -1&1&0&1&-2&2&-1&0 \\ -1&1&0&2&-4&3&-1&0 \\
        -2&2&1&2&-5&4&-2&0 \\ -2&2&1&1&-4&4&-2&0 \\ -1&2&1&0&-3&3&-1&0 \\
        0&2&0&0&-2&2&-1&0 \\ 0&1&0&0&-1&1&0&-1 \end{array}\right)$,
        \\
  $\left(\begin{array}{cccccccc} -1&1&0&1&-2&1&0&0 \\ 0&1&0&1&-3&2&0&0 \\ 0&2&-1&2&-4&2&0&0 \\
        0&2&-1&3&-6&3&1&-1 \\ 0&2&-1&2&-5&3&1&-1 \\ 0&2&-1&1&-3&2&1&-1 \\
        0&1&-1&1&-2&2&0&-1 \\ 0&0&-1&1&-1&1&0&0 \end{array}\right)$,
  $\left(\begin{array}{cccccccc} -1&1&0&1&-2&1&0&0 \\ 0&1&0&1&-3&2&0&0 \\ 0&2&-1&2&-4&2&0&0 \\
        0&3&0&2&-6&3&0&0 \\ 0&3&0&1&-5&3&0&0 \\ 0&3&0&0&-3&2&0&0 \\
        0&2&0&0&-2&2&-1&0 \\ 0&1&0&0&-1&1&0&-1 \end{array}\right)$,
  $\left(\begin{array}{cccccccc} -1&1&0&1&-2&1&0&0 \\ 0&2&-1&1&-2&1&0&0 \\ -1&2&-1&2&-3&2&-1&0 \\
        -1&3&-2&3&-4&2&0&-1 \\ -1&3&-2&2&-3&2&0&-1 \\
        -1&2&-1&1&-2&2&0&-1 \\ 0&1&-1&1&-2&2&0&-1 \\ 0&0&-1&1&-1&1&0&0
        \end{array}\right)$,  \\
  $\left(\begin{array}{cccccccc} 1&-1&0&-1&2&-1&0&0 \\ 0&-2&1&-1&2&-1&0&0 \\
        1&-2&1&-2&3&-2&1&0 \\ 1&-4&1&-2&4&-2&1&0 \\ 1&-4&1&-1&3&-2&1&0 \\
        1&-3&0&0&2&-2&1&0 \\ 0&-2&0&0&2&-2&1&0 \\ 0&-1&0&0&1&-1&0&1 \end{array}\right)$,
   $\left(\begin{array}{cccccccc} 1&-1&0&-1&2&-1&0&0 \\ 0&-2&1&-1&2&-1&0&0 \\
        1&-2&1&-2&3&-2&1&0 \\ 1&-3&2&-3&4&-2&0&1 \\ 1&-3&2&-2&3&-2&0&1 \\
        1&-2&1&-1&2&-2&0&1 \\ 0&-1&1&-1&2&-2&0&1 \\ 0&0&1&-1&1&-1&0&0 \end{array}\right)$,
   $\left(\begin{array}{cccccccc} 1&-1&0&-1&2&-1&0&0 \\ 0&-2&1&-1&2&-1&0&0 \\
        1&-2&1&-2&3&-2&1&0 \\ 1&-3&2&-3&4&-2&1&-1 \\
        1&-3&2&-2&3&-2&1&-1 \\ 1&-2&1&-1&2&-2&1&-1 \\
        0&-1&1&-1&2&-2&1&-1 \\ 0&-1&0&0&1&-1&1&-1
        \end{array}\right)$,  \\
  $\left(\begin{array}{cccccccc} 1&-1&0&-1&2&-1&0&0 \\ 0&-1&0&-1&3&-2&0&0 \\ 0&-2&1&-2&4&-2&0&0 \\
        0&-3&0&-2&6&-3&0&0 \\ 0&-3&0&-1&5&-3&0&0 \\ 0&-3&0&0&3&-2&0&0 \\
        0&-2&0&0&2&-2&1&0 \\ 0&-1&0&0&1&-1&0&1 \end{array}\right)$,
  $\left(\begin{array}{cccccccc} 1&-1&0&-1&2&-1&0&0 \\ 0&-1&0&-1&3&-2&0&0 \\ 0&-2&1&-2&4&-2&0&0 \\
        0&-2&1&-3&6&-3&0&-1 \\ 0&-2&1&-2&5&-3&0&-1 \\
        0&-2&1&-1&3&-2&0&-1 \\ 0&-1&1&-1&2&-2&1&-1 \\
        0&-1&0&0&1&-1&1&-1 \end{array}\right)$,
  $\left(\begin{array}{cccccccc} 1&-1&0&-1&2&-1&0&0 \\ 1&-1&0&-1&2&-2&1&0 \\ 1&-1&0&-2&4&-3&1&0 \\
        2&-1&0&-3&5&-4&2&-1 \\ 2&-1&0&-2&4&-4&2&-1 \\
        1&-1&0&-1&3&-3&1&-1 \\ 0&-1&1&-1&2&-2&1&-1 \\
        0&-1&0&0&1&-1&1&-1 \end{array}\right)$,  \\
  $\left(\begin{array}{cccccccc} 1&-1&0&1&0&-1&0&0 \\ 0&-2&1&2&-1&-1&0&0 \\ 1&-2&1&2&-1&-2&1&0 \\
        1&-4&1&4&-2&-2&1&0 \\ 1&-4&1&3&-1&-2&1&0 \\ 1&-3&0&3&-1&-2&1&0 \\
        0&-2&0&2&0&-2&1&0 \\ 0&-1&0&1&0&-1&0&1 \end{array}\right)$,
  $\left(\begin{array}{cccccccc} 1&-1&0&1&0&-1&0&0 \\ 0&-2&1&2&-1&-1&0&0 \\ 1&-2&1&2&-1&-2&1&0 \\
        1&-3&2&3&-2&-2&0&1 \\ 1&-3&2&2&-1&-2&0&1 \\ 1&-2&1&2&-1&-2&0&1 \\
        0&-1&1&1&0&-2&0&1 \\ 0&0&1&0&0&-1&0&0 \end{array}\right)$,
  $\left(\begin{array}{cccccccc} 1&-1&0&1&0&-1&0&0 \\ 0&-2&1&2&-1&-1&0&0 \\ 1&-2&1&2&-1&-2&1&0 \\
        1&-3&2&3&-2&-2&1&-1 \\ 1&-3&2&2&-1&-2&1&-1 \\
        1&-2&1&2&-1&-2&1&-1 \\ 0&-1&1&1&0&-2&1&-1 \\
        0&-1&0&1&0&-1&1&-1 \end{array}\right)$,  \\
  $\left(\begin{array}{cccccccc} 1&-1&0&1&0&-1&0&0 \\ 0&-1&0&2&0&-2&0&0 \\ 0&-2&1&2&0&-2&0&0 \\
        0&-3&0&4&0&-3&0&0 \\ 0&-3&0&3&1&-3&0&0 \\ 0&-3&0&3&0&-2&0&0 \\
        0&-2&0&2&0&-2&1&0 \\ 0&-1&0&1&0&-1&0&1 \end{array}\right)$,
  $\left(\begin{array}{cccccccc} 1&-1&0&1&0&-1&0&0 \\ 0&-1&0&2&0&-2&0&0 \\ 0&-2&1&2&0&-2&0&0 \\
        0&-2&1&3&0&-3&0&-1 \\ 0&-2&1&2&1&-3&0&-1 \\ 0&-2&1&2&0&-2&0&-1 \\
        0&-1&1&1&0&-2&1&-1 \\ 0&-1&0&1&0&-1&1&-1 \end{array}\right)$,
  $\left(\begin{array}{cccccccc} 1&-1&0&1&0&-1&0&0 \\ 1&-1&0&2&-1&-2&1&0 \\ 1&-1&0&2&0&-3&1&0 \\
        2&-1&0&3&-1&-4&2&-1 \\ 2&-1&0&2&0&-4&2&-1 \\
        1&-1&0&2&0&-3&1&-1 \\ 0&-1&1&1&0&-2&1&-1 \\ 0&-1&0&1&0&-1&1&-1
        \end{array}\right)$,  \\
  $\left(\begin{array}{cccccccc} 1&0&0&0&0&0&0&0 \\ 0&0&1&0&-1&1&0&0 \\ 1&0&1&0&-1&0&1&0 \\
        1&-1&1&1&-2&1&1&0 \\ 1&-1&1&0&-1&1&1&0 \\ 1&0&0&0&-1&1&1&0 \\
        0&0&0&0&0&0&1&0 \\ 0&0&0&0&0&0&0&1 \end{array}\right)$,
  $\left(\begin{array}{cccccccc} 1&0&0&0&0&0&0&0 \\ 0&0&1&0&-1&1&0&0 \\ 1&0&1&0&-1&0&1&0 \\
        1&0&2&0&-2&1&0&1 \\ 1&0&2&-1&-1&1&0&1 \\ 1&1&1&-1&-1&1&0&1 \\
        0&1&1&-1&0&0&0&1 \\ 0&1&1&-1&0&0&0&0 \end{array}\right)$,
  $\left(\begin{array}{cccccccc} 1&0&0&0&0&0&0&0 \\ 0&0&1&0&-1&1&0&0 \\ 1&0&1&0&-1&0&1&0 \\
        1&0&2&0&-2&1&1&-1 \\ 1&0&2&-1&-1&1&1&-1 \\ 1&1&1&-1&-1&1&1&-1 \\
        0&1&1&-1&0&0&1&-1 \\ 0&0&0&0&0&0&1&-1 \end{array}\right)$,
        \\
  $\left(\begin{array}{cccccccc} -1&-1&1&0&1&-2&1&0 \\ -1&-2&0&1&1&-2&1&0 \\ -2&-2&1&1&1&-3&1&1 \\
        -3&-4&1&2&2&-4&1&1 \\ -3&-3&1&1&2&-3&1&1 \\ -2&-2&1&0&2&-2&0&1 \\
        -2&-1&1&0&1&-1&0&0 \\ -1&-1&1&0&0&0&0&0 \end{array}\right)$,
  $\left(\begin{array}{cccccccc} 1&0&0&0&0&0&0&0 \\ 0&1&0&0&0&0&0&0 \\ 0&0&1&0&0&0&0&0 \\
        0&1&1&0&0&0&0&-1 \\ 0&1&1&-1&1&0&0&-1 \\ 0&1&1&-1&0&1&0&-1 \\
        0&1&1&-1&0&0&1&-1 \\ 0&0&0&0&0&0&1&-1 \end{array}\right)$,
  $\left(\begin{array}{cccccccc} 1&0&0&0&0&0&0&0 \\ 1&1&0&0&-1&0&1&0 \\ 1&1&0&0&0&-1&1&0 \\
        2&2&0&0&-1&-1&2&-1 \\ 2&2&0&-1&0&-1&2&-1 \\ 1&2&0&-1&0&0&1&-1 \\
        0&1&1&-1&0&0&1&-1 \\ 0&0&0&0&0&0&1&-1 \end{array}\right)$.
%%\end {document}

The character table of $G^{s_{96}}$:\\
\noindent \begin{tabular}{c|cccccccccccccccccccccccc}
  & & & & & & & & & &10 & & & & & & & & & & 20& & & &\\\hline
$\chi_{96}^{(1)}$&1&1&1&1&1&1&1&1&1&1&1&1&1&1&1&1 &1&1&1&1&1&1&1&1

\\  $\chi_{96}^{(2)}$&1&-1&-1&-1&-1&-1&-1&-1&-1&-1&-1&-1&1&1&1&1
&1&1&1&1&1&1&1&1

\\  $\chi_{96}^{(3)}$&1&1&1&-1&-1&-1&-1&-1&-1&-1&-1&-1&-1&-1&-1&1
&1&1&1&1&1&1&1&1

\\  $\chi_{96}^{(4)}$&1&1&1&1&1&1&1&1&1&1&1&1&-1&-1&-1&-1
&-1&-1&-1&-1&-1&-1&-1&-1

\\  $\chi_{96}^{(5)}$&1&-1&-1&1&1&1&1&1&1&1&1&1&1&1&1&-1
&-1&-1&-1&-1&-1&-1&-1&-1

\\  $\chi_{96}^{(6)}$&1&1&1&-1&-1&-1&-1&-1&-1&-1&-1&-1&1&1&1&-1
&-1&-1&-1&-1&-1&-1&-1&-1

\\  $\chi_{96}^{(7)}$&1&-1&-1&-1&-1&-1&-1&-1&-1&-1&-1&-1&-1&-1&-1&-1
&-1&-1&-1&-1&-1&-1&-1&-1

\\
\end {tabular}

\noindent \begin{tabular}{c|cccccccccccccccccccccccc}
  & & & & & & & & & &10 & & & & & & & & & & 20& & & &\\\hline

\\  $\chi_{96}^{({46})}$&2&.&.&-B&B&.&.&-B&B&B&.&-B&.&.&.&B
&-B&.&.&B&-B&-B&.&B

\\  $\chi_{96}^{({47})}$&2&.&.&C&B&.&.&-/C&-C&/C&.&-B&.&.&.&-C
&-B&.&.&/C&C&-/C&.&B

\\  $\chi_{96}^{({48})}$&2&.&.&/C&-B&.&.&-C&-/C&C&.&B&.&.&.&-/C
&B&.&.&C&/C&-C&.&-B

\\  $\chi_{96}^{({49})}$&2&.&.&-C&-B&.&.&/C&C&-/C&.&B&.&.&.&C
&B&.&.&-/C&-C&/C&.&-B

\\  $\chi_{96}^{({50})}$&2&.&.&-/C&B&.&.&C&/C&-C&.&-B&.&.&.&/C
&-B&.&.&-C&-/C&C&.&B

\\  $\chi_{96}^{({51})}$&2&.&.&/A&-1&D&2&A&/A&A&/D&-1&.&.&.&/A
&-1&D&2&A&/A&A&/D&-1

\\
\end{tabular}

\noindent \begin{tabular}{c|cccccccccccccccccccccccc}
  & & & & & &30 & & & & & & & & & &40 & & & & & & & &\\\hline
$\chi_{96}^{(1)}$&1&1&1&1&1&1&1&1&1&1&1&1&1&1&1&1 &1&1&1&1&1&1&1&1

\\  $\chi_{96}^{(2)}$&1&1&1&1&1&1&1&1&1&-1&-1&-1&-1&-1&-1&-1
&-1&-1&-1&-1&-1&1&1&1

\\  $\chi_{96}^{(3)}$&-1&-1&-1&1&1&1&1&1&1&1&1&1&-1&-1&-1&-1
&-1&-1&1&1&1&-1&-1&-1

\\  $\chi_{96}^{(4)}$&1&1&1&1&1&1&1&1&1&-1&-1&-1&-1&-1&-1&-1
&-1&-1&-1&-1&-1&1&1&1

\\  $\chi_{96}^{(5)}$&-1&-1&-1&1&1&1&1&1&1&1&1&1&-1&-1&-1&-1
&-1&-1&1&1&1&-1&-1&-1

\\  $\chi_{96}^{(6)}$&-1&-1&-1&1&1&1&1&1&1&-1&-1&-1&1&1&1&1
&1&1&-1&-1&-1&-1&-1&-1

\\  $\chi_{96}^{(7)}$&1&1&1&1&1&1&1&1&1&1&1&1&1&1&1&1
&1&1&1&1&1&1&1&1

\\  $\chi_{96}^{(8)}$&-1&-1&-1&1&1&1&1&1&1&-1&-1&-1&1&1&1&1
&1&1&-1&-1&-1&-1&-1&-1

\\  $\chi_{96}^{(9)}$&-A&-/A&1&-A&1&-/A&-A&1&-/A&A&/A&-1&A&-1&/A&A
&-1&/A&A&/A&-1&-A&-/A&1

\\  $\chi_{96}^{({10})}$&-/A&-A&1&-/A&1&-A&-/A&1&-A&/A&A&-1&/A&-1&A&/A
&-1&A&/A&A&-1&-/A&-A&1

\\  $\chi_{96}^{({11})}$&A&/A&-1&-A&1&-/A&-A&1&-/A&-A&-/A&1&A&-1&/A&A
&-1&/A&-A&-/A&1&A&/A&-1

\\  $\chi_{96}^{({12})}$&/A&A&-1&-/A&1&-A&-/A&1&-A&-/A&-A&1&/A&-1&A&/A
&-1&A&-/A&-A&1&/A&A&-1

\\  $\chi_{96}^{({13})}$&-/A&-A&1&-/A&1&-A&-/A&1&-A&/A&A&-1&/A&-1&A&/A
&-1&A&/A&A&-1&-/A&-A&1

\\  $\chi_{96}^{({14})}$&-A&-/A&1&-A&1&-/A&-A&1&-/A&A&/A&-1&A&-1&/A&A
&-1&/A&A&/A&-1&-A&-/A&1

\\  $\chi_{96}^{({15})}$&/A&A&-1&-/A&1&-A&-/A&1&-A&-/A&-A&1&/A&-1&A&/A
&-1&A&-/A&-A&1&/A&A&-1
\end{tabular}

\noindent \begin{tabular}{c|cccccccccccccccccccccccc}
  & & & & & &30 & & & & & & & & & &40 & & & & & & & &\\\hline
\\  $\chi_{96}^{({16})}$&A&/A&-1&-A&1&-/A&-A&1&-/A&-A&-/A&1&A&-1&/A&A
&-1&/A&-A&-/A&1&A&/A&-1

\\  $\chi_{96}^{({17})}$&A&/A&-1&-A&1&-/A&-A&1&-/A&A&/A&-1&-A&1&-/A&-A
&1&-/A&A&/A&-1&A&/A&-1

\\  $\chi_{96}^{({18})}$&/A&A&-1&-/A&1&-A&-/A&1&-A&/A&A&-1&-/A&1&-A&-/A
&1&-A&/A&A&-1&/A&A&-1

\\  $\chi_{96}^{({19})}$&-A&-/A&1&-A&1&-/A&-A&1&-/A&-A&-/A&1&-A&1&-/A&-A
&1&-/A&-A&-/A&1&-A&-/A&1

\\  $\chi_{96}^{({20})}$&-/A&-A&1&-/A&1&-A&-/A&1&-A&-/A&-A&1&-/A&1&-A&-/A
&1&-A&-/A&-A&1&-/A&-A&1

\\  $\chi_{96}^{({21})}$&/A&A&-1&-/A&1&-A&-/A&1&-A&/A&A&-1&-/A&1&-A&-/A
&1&-A&/A&A&-1&/A&A&-1

\\
\end{tabular}%

\noindent where A = -E(3)$^2$
  = (1+ER(-3))/2 = 1+b3,
B = -E(3)+E(3)$^2$
  = -ER(-3) = -i3,
C = E(3)+2*E(3)$^2$
  = (-3-ER(-3))/2 = -2-b3,
D = 2*E(3)$^2$
  = -1-ER(-3) = -1-i3,
E = 4*E(3)$^2$
  = -2-2*ER(-3) = -2-2i3.

The generators of $G^{s_{97}}$ are:\\
$\left(\begin{array}{cccccccc} -1&1&1&-2& 1&0&1&0 \\
0&1&1&-3&1&1&1&0 \\ 0&2&1&-4&2&0&2&0
\\ 0&3&2&-6&2&1&3&-1 \\ 0&2&2&-5&2&0&3&-1 \\ 0&
2&1&-4&2&0&2&-1 \\
        0&1&1&-3&2&0&1&-1 \\ 0&0&0&-1&1&0&1&-1 \end{array}\right)$,
  $\left(\begin{array}{cccccccc} -1&0&1&0&0&0&0&0 \\ 0&1&0&0&0&0&0&0 \\ 0&0&1&0&0&0&0&0 \\
        0&0&0&1&0&0&0&0 \\ 0&0&0&0&1&0&0&0 \\ 0&0&0&0&0&1&0&0 \\
        0&0&0&0&0&0&1&0 \\ 0&0&0&0&0&0&0&1 \end{array}%
\right)$.

The representatives of conjugacy classes of   $G^{s_{97}}$ are:\\
$\left(%
\right)$.
%%\end {document}

The character table of $G^{s_{97}}$:\\
\noindent %

\noindent where A = E(3)$^2$
  = (-1-ER(-3))/2 = -1-b3,
B = 2*E(3)$^2$
  = -1-ER(-3) = -1-i3,
C = 3*E(3)$^2$
  = (-3-3*ER(-3))/2 = -3-3b3,
D = -E(4)
  = -ER(-1) = -i,
E = -E(12)$^{11}$,F = -2*E(4)
  = -2*ER(-1) = -2i,
G = -2*E(12)$^{11}$,H = -3*E(4)
  = -3*ER(-1) = -3i,
I = -3*E(12)$^{11}$.

The generators of $G^{s_{98}}$ are:\\
$\left(%
\right).$
%%\end {document}

The character table of $G^{s_{98}}$:\\
\noindent %

\noindent where A = -E(4)
  = -ER(-1) = -i,
B = -E(12)$^7$,C = -E(3)
  = (1-ER(-3))/2 = -b3.

The generators of $G^{s_{99}}$ are:\\
$\left(%
\right)$.

The character table of $G^{s_{99}}$:\\
${}_{
 \!\!\!\!\!_{
\!\!\!\!\!_{ \!\!\!\!\!  _{  \small
 \noindent
%

\noindent where A = E(4)
  = ER(-1) = i,
B = -E(3)
  = (1-ER(-3))/2 = -b3,
C = E(12)$^7$.

The generators of $G^{s_{100}}$ are:\\
$\left(%
\right)$.
%\end {document}

The character table of $G^{s_{100}}$:\\
\noindent %

\noindent where A = E(3)$^2$
  = (-1-ER(-3))/2 = -1-b3,
B = 2*E(3)$^2$
  = -1-ER(-3) = -1-i3,
C = 3*E(3)$^2$
  = (-3-3*ER(-3))/2 = -3-3b3.

The generators of $G^{s_{101}}$ are:\\
$\left(%
\right)$.

The character table of $G^{s_{101}}$:\\
${}_{
 \!\!\!\!\!_{
\!\!\!\!\!_{ \!\!\!\!\!  _{  \small
 \noindent
%

\noindent where A = -E(4)
  = -ER(-1) = -i,
B = -E(3)$^2$
  = (1+ER(-3))/2 = 1+b3,
C = E(12)$^{11}$.

The generators of $G^{s_{102}}$ are:\\
$\left(%
\right)$.
%%\end {document}

The character table of $G^{s_{102}}$:\\
\noindent %

\noindent where A = E(3)$^2$
  = (-1-ER(-3))/2 = -1-b3,
B = 2*E(3)
  = -1+ER(-3) = 2b3,
C = 4*E(3)
  = -2+2*ER(-3) = 4b3.

The generators of $G^{s_{103}}$ are:\\
$\left(%
\right)$.

The character table of $G^{s_{103}}$:\\
${}_{
 \!\!\!\!\!\!_{
\!\!\!\!\!\!_{ \!\!\!\!\!\!  _{  \small
 \noindent
%

\noindent where A = -E(3)$^2$
  = (1+ER(-3))/2 = 1+b3,
B = E(4)
  = ER(-1) = i,
C = E(12)$^{11}$,D = -2*E(4)
  = -2*ER(-1) = -2i,
E = 3*E(4)
  = 3*ER(-1) = 3i,
F = 2*E(3)
  = -1+ER(-3) = 2b3,
G = 3*E(3)
  = (-3+3*ER(-3))/2 = 3b3,
H = -2*E(12)$^{11}$,I = -3*E(12)$^{11}$.

The generators of $G^{s_{104}}$ are:\\
$\left(%
\right)$.
%%\end {document}

The character table of $G^{s_{104}}$:\\
${}_{
 \!\!\!\!\!\!_{
\!\!\!\!\!\!_{ \!\!\!\!\!\!  _{  \small
 \noindent
%

\noindent where A = -E(4)
  = -ER(-1) = -i,
B = E(8),C = -2*E(4)
  = -2*ER(-1) = -2i,
D = 2*E(8)$^3$.

The generators of $G^{s_{105}}$ are:\\
$\left(%
\right)$.

%%\end {document}

The character table of $G^{s_{105}}$:\\
${}_{
 \!\!\!\!_{
\!\!\!\!_{ \!\!\!\!  _{  \small
 \noindent
%

\noindent where A = -E(5)$^3$,B = -E(5),C = E(3)
  = (-1+ER(-3))/2 = b3,
D = E(15)$^7$,E = E(15)$^{14}$,F = E(15)$^{13}$,G = E(15)$^4$.

The generators of $G^{s_{106}}$ are:\\
$\left(%
\right).$

The representatives of conjugacy classes of   $G^{s_{106}}$ are:\\
$\left(%
\right).$

The character table of $G^{s_{106}}$:\\
\noindent %

\noindent where A = E(3)$^2$
  = (-1-ER(-3))/2 = -1-b3,
B = 5*E(3)$^2$
  = (-5-5*ER(-3))/2 = -5-5b3,
C = 9*E(3)$^2$
  = (-9-9*ER(-3))/2 = -9-9b3,
D = 10*E(3)$^2$
  = -5-5*ER(-3) = -5-5i3,
E = 16*E(3)$^2$
  = -8-8*ER(-3) = -8-8i3,
F = -3*E(3)$^2$
  = (3+3*ER(-3))/2 = 3+3b3,
G = 2*E(3)$^2$
  = -1-ER(-3) = -1-i3.
%%\end {document}

The generators of $G^{s_{107}}$ are:\\
$\left(%
\right).$

The representatives of conjugacy classes of   $G^{s_{107}}$ are:\\
$\left(%
\right).$

The character table of $G^{s_{107}}$:\\
\noindent %

\noindent where A = E(8)$^3$,B = -E(4)
  = -ER(-1) = -i.

The generators of $G^{s_{108}}$ are:\\
$\left(%
\right)$.
%%\end {document}
The character table of $G^{s_{108}}$:\\
\noindent %

\noindent where A = E(3)$^2$
  = (-1-ER(-3))/2 = -1-b3,
B = 2*E(3)$^2$
  = -1-ER(-3) = -1-i3,
C = 4*E(3)$^2$
  = -2-2*ER(-3) = -2-2i3,
D = 6*E(3)$^2$
  = -3-3*ER(-3) = -3-3i3,
E = 8*E(3)$^2$
  = -4-4*ER(-3) = -4-4i3,
F = 9*E(3)$^2$
  = (-9-9*ER(-3))/2 = -9-9b3,
G = 12*E(3)$^2$
  = -6-6*ER(-3) = -6-6i3,
H = 16*E(3)$^2$
  = -8-8*ER(-3) = -8-8i3,
I = -3*E(3)$^2$
  = (3+3*ER(-3))/2 = 3+3b3.

The generators of $G^{s_{109}}$ are:\\
$\left(%
\right).$

The character table of $G^{s_{109}}$:\\
\noindent %

\noindent where A = E(3)$^2$
  = (-1-ER(-3))/2 = -1-b3,
B = 2*E(3)$^2$
  = -1-ER(-3) = -1-i3,
C = 3*E(3)$^2$
  = (-3-3*ER(-3))/2 = -3-3b3.

The generators of $G^{s_{110}}$ are:\\
$\left(%
\right).$

The representatives of conjugacy classes of   $G^{s_{110}}$ are:\\
$\left(%
\right).$

The character table of $G^{s_{110}}$:\\
${}_{
 \!\!\!\!\!_{
\!\!\!\!\!_{ \!\!\!\!\!  _{  \small
 \noindent
%
 }}}}$

\noindent where A = -E(4)
  = -ER(-1) = -i,
B = -E(12)$^7$,C = E(3)$^2$
  = (-1-ER(-3))/2 = -1-b3,
D = E(8),E = E(24)${11}$,F = E(24)$^{19}$.

The generators of $G^{s_{111}}$ are:\\
$\left(%
\right).$

The character table of $G^{s_{111}}$:\\
\noindent %

\noindent where A = E(3)$^2$
  = (-1-ER(-3))/2 = -1-b3,
B = E(4)
  = ER(-1) = i,
C = -2*E(4)
  = -2*ER(-1) = -2i,
D = -E(12)$^7$,E = -2*E(3)
  = 1-ER(-3) = 1-i3,
F = -2*E(12)$^7$.

The generators of $G^{s_{112}}$ are:\\
$\left(%
\right).$
%%\end {document}
The character table of $G^{s_{112}}$:\\
\noindent %

\noindent where A = E(3)
  = (-1+ER(-3))/2 = b3,
B = -2*E(3)$^2$
  = 1+ER(-3) = 1+i3.

\vskip 0.3cm

\noindent{\large\bf Acknowledgement}: We would like to thank Prof.
N. Andruskiewitsch and Dr. F. Fantino for suggestions and help.

}
\end {document}